
\input amstex
\documentstyle{amsppt}

\magnification=\magstep1
\hsize=6.5truein
\vsize=9truein

\document

\baselineskip=12pt

\font \smallit=cmti10 at 9pt
\font \smallsl=cmsl10 at 9pt
\font \bfit=cmmib9 at 10pt

\def \loongrightarrow {\relbar\joinrel\relbar\joinrel\rightarrow}
\def \llongrightarrow
{\relbar\joinrel\relbar\joinrel\relbar\joinrel\rightarrow}

%

%
\def \gerg {\frak g}

\def \gerh {\frak h}

%
\def \N {\Bbb N}

\def \HA {{\Cal H}{\Cal A}}

\def \Fstel {{F\!\phantom{)}}^{\scriptscriptstyle \!\circledast}}
\def \Finfty {{F\!\phantom{)}}^{{\scriptscriptstyle \!\infty}}}
\def \Fhstel {{F_h\!\!\!\phantom{)}}^{\scriptscriptstyle \!\circledast}}
\def \Fhinfty {{F_h\!\!\!\phantom{)}}^{{\scriptscriptstyle \!\infty}}}
\def \otimeshat {\,\widehat{\otimes}\,}
\def \otimestilde {\,\widetilde{\otimes}\,}
\def \otimesbar {\,\bar{\otimes}\,}
\def \otimescheck {\,\check{\otimes}\,}
\def \QUEA {\hbox{$ \displaystyle{\Cal{Q\hskip1ptUE\hskip-1ptA}} $}}
\def \QFSHA {\hbox{$ \displaystyle{\Cal{Q\hskip1ptFSHA}} $}}
\def \lsF {\hbox{$ {\phantom{\big|}}^{\scriptscriptstyle F\!\!} $}}


\topmatter

{\ }

\vskip-33pt

\hfill   {{\sl Annales de l'Institut Fourier\/}  {\bf 52},
no.~3 (2002), 809--834}
\hskip19pt   {\ }

\vskip3pt

\hfill   {\tt http://aif.cedram.org/item?id=AIF{\underscore}2002{\underscore}{\underscore}52{\underscore}3{\underscore}809{\underscore}0}
\hskip19pt   {\ }

\vskip31pt

\title
  The quantum duality principle
\endtitle

\author
       Fabio Gavarini
\endauthor


\affil
  Universit\`a degli Studi di Roma ``Tor Vergata'' ---
Dipartimento di Matematica  \\
  Via della Ricerca Scientifica 1, I-00133 Roma --- ITALY  \\
\endaffil

\address\hskip-\parindent
   Universit\`a degli Studi di Roma ``Tor Vergata''
---   Dipartimento di Matematica  \newline
   Via della Ricerca Scientifica 1, I-00133 Roma, ITALY
---   e-mail: \ gavarini\@{}mat.uniroma2.it
\endaddress

\abstract
  The "quantum duality principle" states that the quantisation of a Lie
bialgebra   --- via a quantum universal enveloping algebra (in short,
QUEA) ---   also provides a quantisation of the dual Lie bialgebra
(through its associated formal Poisson group)   --- via a quantum
formal series Hopf algebra (QFSHA) ---   and, conversely, a QFSHA
associated to a Lie bialgebra (via its associated formal Poisson
group) yields a QUEA for the dual Lie bialgebra as well; more in
detail, there exist functors  $ \; \QUEA \longrightarrow \QFSHA\; $
and  $ \; \QFSHA \longrightarrow \QUEA \, $,  \, inverse to each other,
such that in both cases the Lie bialgebra associated to the target
object is the dual of that of the source object.  Such a result was
claimed true by Drinfeld, but seems to be unproved in the literature:
I give here a thorough detailed proof of it.
\endabstract

\endtopmatter

\footnote""{Keywords: \ {\sl Quantum Groups, topological Hopf
algebras}.}

\footnote""{ 2000 {\it Mathematics Subject Classification:}
\  Primary 17B37, 20G42, Secondary 81R50, 16W30. }

%
%

\hfill  \hbox{\vbox{ \baselineskip=10pt
          \hbox{\smallit  \   "Dualitas dualitatum }
          \hbox{\smallit \ \ \;\, et omnia dualitas" }
             \vskip4pt
          \hbox{\smallsl    N.~Barbecue, "Scholia" } } \hskip1truecm }

\vskip5pt

\centerline {\bf  Introduction }

\vskip10pt

   The  {\it quantum duality principle\/}  is known in literature
under at least two formulations.  One claims that quantum function
algebras associated to dual Poisson groups can be considered to be
dual   --- in the Hopf sense ---   to each other; and similarly for
quantum enveloping algebras (cf.~[FRT] and [Se]).  The second one, due
to Drinfeld (cf.~[Dr]), states that any quantisation of the universal
enveloping algebra of a Poisson group can also be understood   --- in
some sense ---   as a quantisation of the dual formal Poisson group,
and, conversely, any quantisation of a formal Poisson group also
"serves" as a quantisation of the universal enveloping algebra of
the dual Poisson group: this is the point of view we are interested
in.  I am now going to describe this result more in detail.
                                            \par
   Let  $ \Bbbk $  be a field of zero characteristic.  Let  $ \gerg $
be a finite dimensional Lie algebra over  $ \Bbbk $,  $ U(\gerg) $  its
universal enveloping algebra: then  $ U (\gerg) $  has a natural
structure of Hopf algebra.  Let  $ F[[\gerg]] $  be the (algebra
of regular functions on the) formal group associated to  $ \gerg $:
it is a complete topological Hopf algebra (the coproduct taking values
in a suitable topological tensor product of the algebra with itself),
which has two realisations.  The first one is as follows: if  $ G $
is an affine algebraic group with tangent Lie algebra  $ \gerg $,
and  $ F[G] $  is the algebra of regular functions on  $ G $,  then
$ F[[\gerg]] $  is the  $ {\frak m}_e $--completion  of  $ F[G] $  at
the maximal ideal  $ {\frak m}_e $  of the identity element  $ \, e \in
G \, $,  endowed with the  $ {\frak m}_e $--adic  topology.  The second
one is  $ \, F[[\gerg]] := {U(\gerg)}^* \, $,  \, the linear dual of
$ U(\gerg) \, $),  endowed with the weak topology.  In any case,
$ U(\gerg) $  identifies with the  {\sl topological dual\/}  of
$ F[[\gerg]] $,  i.e.~the set of all  $ \Bbbk $--linear  continuous maps
from  $ F[[\gerg]] $  to  $ \Bbbk $,  where  $ \Bbbk $  is given the discrete
topology; similarly  $ \, F[[\gerg]] = {U(\gerg)}^* \, $  is also the
topological dual of  $ U(\gerg) $  if we take on the latter space the
discrete topology: in particular, a (continuous) biduality theorem
relates  $ U(\gerg ) $  and $ F[[\gerg]] $,  and evaluation yields a
natural Hopf pairing among them.  Now assume  $ \gerg $  is a  {\sl
Lie bialgebra\/}: then  $ U(\gerg) $  is a  {\sl co-Poisson Hopf\/}
algebra, $ F[[\gerg]] $  is a topological  {\sl Poisson\/}  Hopf
algebra, and the above pairing is compatible with these additional
co-Poisson and Poisson structures.  Further, the dual  $ \gerg^* $
of  $ \gerg $ is a Lie bialgebra as well, so we can consider also
$ U(\gerg^*) $  and  $ F[[\gerg^*]] $.
                                            \par
   Let  $ \gerg $  be a Lie bialgebra.  A quantisation of  $ U(\gerg) $
is, roughly speaking, a topological Hopf  $ \Bbbk[[h]] $--algebra  which
for  $ \, h = 0 \, $  is isomorphic, as a co-Poisson Hopf algebra,
to  $ U(\gerg) $:  these objects form a category, called \QUEA.
Similarly, a quantisation of  $ F[[\gerg]] $  is, in short, a
topological Hopf  $ \Bbbk[[h]] $--algebra  which for  $ \, h = 0
\, $  is isomorphic, as a topological Poisson Hopf algebra, to
$ F[[\gerg]] $:  we call \QFSHA{} the category formed by these objects.
                                            \par
   The quantum duality principle (after Drinfeld) states that
there exist two functors, namely  $ \; {(\ )}' \colon \, \QUEA
\longrightarrow \QFSHA \; $  and  $ \; {(\ )}^\vee \colon \, \QFSHA
\longrightarrow \QUEA \, $,  \, which are inverse of each other, and
if  $ \, U_h(\gerg) \, $  is a quantisation of  $ \, U(\gerg) \, $
and  $ \, F_h[[\gerg]] \, $  is a quantisation of  $ \, F[[\gerg]]
\, $,  then  $ \, {U_h(\gerg)}' \, $  is a quantisation of  $ \,
F[[\gerg^*]] \, $,  and  $ \, {F_h[[\gerg]]}^\vee \, $  is a
quantisation of  $ \, U(\gerg^*) \, $.
                                            \par
   This paper provides an explicit thorough proof (seemingly, the
first one in the literature) of this result.  I also point out
some further details and what is true when  $ \Bbbk $  has positive
characteristic, and sketch a plan for generalizing all this to
the infinite dimensional case.
                                            \par
   Note that several properties of the objects I consider have been
discovered and exploited in the works by Etingof and Kazhdan (see [EK1],
[EK2]), by Enriquez (cf.~[E]) and by Kassel and Turaev (cf.~[KT]), who
deal with some  {\sl special cases\/}  of quantum groups, arising from
a specific construction, and also applied Drinfeld's results.  The
analysis in the present paper shows that that those properties are
often direct consequences of more general facts.
                                            \par
   I point out that Drinfeld's result is essentially  {\sl local\/}
in nature, as it deals with quantisations over the ring of formal
series and ends up only with infinitesimal data, i.e.~objects attached
to Lie bialgebras; a {\sl global\/}  version of the principle, dealing
with quantum groups over a ring of Laurent polynomials, which give
information on the global data of the underlying Poisson groups will
be provided in a forthcoming paper (cf.~[Ga2]): this is useful in
applications, e.g.~it yields a quantum duality
      principle for Poisson homogeneous spaces, cf.~[CG]. \break

\vskip4pt

\centerline{ ACKNOWLEDGEMENTS }

\vskip2pt

  \hbox{The author thanks Pierre Baumann and
Alessandro D'Andrea for their valuable remarks.}

\vskip1,3truecm

\centerline {\bf \S\; 1 \ Notation and terminology }

\vskip10pt

  {\bf 1.1 Topological  $ \Bbbk[[h]] $--modules  and topological Hopf
$ \Bbbk[[h]] $--algebras.} \, Let  $ \Bbbk $  be a fixed field,  $ h $  an
indeterminate.  The ring  $ \Bbbk[[h]] $  will always be considered as
a topological ring w.r.t.~the  $ h $--adic  topology.  Let  $ X $
be any  $ \Bbbk[[h]] $--module.  We set  $ \, X_0 := X \big/ h X = \Bbbk
\otimes_{\Bbbk[[h]]} X \, $, \, a  $ \Bbbk $--module  (via scalar restriction
$ \, \Bbbk[[h]] \rightarrow \Bbbk[[h]] \big/ h \, \Bbbk[[h]] \cong \Bbbk \, $)  which
we call the  {\sl specialisation\/}  of  $ X $  at  $ \, h = 0 \, $,
or  {\sl semiclassical limit\/}  of  $ X \, $;  we shall also use
notation  $ \, X \,{\buildrel \, h \rightarrow 0 \, \over
\llongrightarrow}\, \overline{Y} \, $  to mean  $ \, X_0 \cong
\overline{Y} \, $.  Note that if  $ X $  is a topological
$ \Bbbk[[h]] $--module  which is torsionless, complete and separated
w.r.t.~the  $ h $--adic  topology then there is a natural isomorphism
of  $ \Bbbk[[h]] $--modules  $ \, X \cong X_0[[h]] \, $:  \, indeed,
choose any  $ \Bbbk $--basis  $ {\big\{ b_i \big\}}_{i \in I} $  of
$ X_0 $,  and pick any subset  $ {\big\{ \beta_i \big\}}_{i \in I}
\subseteq X \, $  such that  $ \, \beta_i \mod h = b_i \; (\forall\,
i \,) \, $:  \, then an isomorphism as required is given by  $ \,
\beta_i \mapsto b_i \, $  (however, topologies on either side may
be different).
                                            \par
   For later use, we also set  $ \, \lsF X := \Bbbk((h))
\otimes_{\Bbbk[[h]]} X \, $,  \, a vector space over
$ \Bbbk((h)) $,  which is not equipped with any topology.
%
%
%
%
                                            \par
   If  $ X $  is a topological  $ \Bbbk[[h]] $--module,  we let its
{\it full dual\/}  to be  $ \, X^* := \text{\it Hom\,}_{\Bbbk[[h]]}
\big( X, \Bbbk[[h]] \big) \, $,  \, and its  {\it topological dual\/}
to be  $ \, X^\star := \big\{\, f \in X^* \,\big\vert\, f \text{\
is continuous} \,\big\} \, $.  Note that  $ \, X^* = X^\star \, $
when the topology on  $ X $  is the  $ h $--adic  one.
                                            \par
   We introduce now two tensor categories of topological
$ \Bbbk[[h]] $--modules,  $ {\Cal T}_{\otimeshat} $  and
$ {\Cal P}_{\otimestilde} $:  the first one is modeled on the tensor
category of  {\sl discrete\/}  topological  $ \Bbbk $--vector
spaces, the second one is modeled on the category of  {\sl
linearly compact\/}  topological  $ \Bbbk $--vector  spaces.
                                            \par
   Let  $ {\Cal T}_{\otimeshat} $  be the category whose objects are
all topological  $ \Bbbk[[h]] $--modules  which are topologically free
(i.e.~isomorphic to  $ V[[h]] $  for some  $ \Bbbk $--vector  space
$ V $,  with the  $ h $--adic  topology) and whose morphisms are
the  $ \Bbbk[[h]] $--linear  maps (which are automatically continuous).
This is a tensor category w.r.t.~the tensor product  $ \, T_1 \otimeshat
T_2 \, $  defined to be the separated  $ h $--adic  completion of the
algebraic tensor product  $ \, T_1 \otimes_{\Bbbk[[h]]} T_2 \, $  (for
all  $ T_1 $,  $ T_2 \in {\Cal T}_{\otimeshat} $).
                                            \par
   Let  $ {\Cal P}_{\otimestilde} $  be the category whose objects
are all topological  $ \Bbbk[[h]] $--modules  isomorphic to modules
of the type  $ {\Bbbk[[h]]}^E $  (the Cartesian product indexed by
$ E $,  with the Tikhonov product topology) for some set  $ E \, $:
\, these are complete w.r.t.~to the weak topology, in fact they are
isomorphic to the projective limit of their finite free submodules
(each one taken with the  $ h $--adic  topology); the morphisms in
$ {\Cal P}_{\otimestilde} $  are the  $ \Bbbk[[h]] $--linear  continuous
maps.  This is a tensor category w.r.t.~the tensor product  $ \, P_1
\otimestilde P_2 \, $  defined to be the completion of the algebraic
tensor product  $ \, P_1 \otimes_{\Bbbk[[h]]} P_2 \, $  w.r.t.~the weak
topology:  therefore  $ \, P_i \cong {\Bbbk[[h]]}^{E_i} $  ($ i = 1 $,
$ 2 $)  yields  $ \, \, P_1 \otimestilde P_2 \cong {\Bbbk[[h]]}^{E_1
\times E_2} \, $  (for all $ P_1 $,  $ P_2 \in
{\Cal P}_{\otimestilde} $).
                                            \par
   Note that the objects of  $ {\Cal T}_{\otimeshat} $  and of
$ {\Cal P}_{\otimestilde} $  are complete and separated w.r.t.~the
$ h $--adic  topology, so by the previous remark one has  $ \, X
\cong X_0[[h]] \, $  for each of them.
                                            \par
   We denote by  $ \HA_{\otimeshat} $  the subcategory of
$ {\Cal T}_{\otimeshat} $  whose objects are all the Hopf
algebras in  $ {\Cal T}_{\otimeshat} $  and whose morphisms
are all the Hopf algebra morphisms in $ {\Cal T}_{\otimeshat} $.
Similarly, we call  $ \HA_{\otimestilde} $  the subcategory
of  $ {\Cal P}_{\otimestilde} $  whose objects are all the Hopf
algebras in  $ {\Cal P}_{\otimestilde} $  and whose morphisms
are all the Hopf algebra morphisms in  $ {\Cal P}_{\otimestilde} $.
Moreover, we define  $ \, \HA_{\otimestilde}^{w-I} \, $  to be the
full subcategory of  $ \, \HA_{\otimestilde} $  whose objects are
all the  $ \, H \in \HA_{\otimestilde} \, $  whose (weak) topology
coincides with the  $ I_{\scriptscriptstyle H} $--adic  topology,
where  $ I_{\scriptscriptstyle H} := h \, H + \text{\it Ker}\,
(\epsilon) = \epsilon^{-1} \big( h \, \Bbbk[[h]] \big) \, $.
                                            \par
   As a matter of notation, when dealing with a (possibly topological)
Hopf algebra  $ H $,  I shall denote by  $ m $  its product, by
$ 1 $  its unit element, by  $ \Delta $  its coproduct, by
$ \epsilon $  its counit and by  $ S $  its antipode; subscripts
$ H $  will be added whenever needed for clarity.  Note that the
objects of  $ \HA_{\otimeshat} $  and of  $ \HA_{\otimestilde} $
are  {\sl topological\/}  Hopf algebras, not standard ones:
in particular, in  $ \sigma $--notation  $ \, \Delta(x) =
\sum_{(x)} x_{(1)} \otimes x_{(2)} \, $  the sum is
\hbox{understood in topological sense.}

\vskip7pt

\proclaim{Definition 1.2}  (cf.~[Dr], \S~7)
                                         \hfill\break
  \indent  (a) \, We call  {\sl quantized universal enveloping
algebra\/}  (in short, QUEA)  any  $ \, H \in \HA_{\otimeshat} $
such that  $ \, H_0 := H \big/ h H \, $  is a co-Poisson Hopf algebra
isomorphic to  $ U(\gerg) $  for some finite dimensional Lie bialgebra
$ \gerg $  (over  $ \Bbbk $);  in this case we write  $ \, H =
U_h(\gerg) \, $,  \, and say  $ H $  is a  {\sl quantisation\/}
of  $ U(\gerg) $.  We call \QUEA{} the full subcategory of
$ \HA_{\otimeshat} $  whose objects are QUEA, relative to
all possible  $ \gerg $  (see also  Remark 1.3{\it (a)\/}  below).
                                         \hfill\break
  \indent  (b) \, We call  {\sl quantized formal series Hopf algebra\/}
(in short, QFSHA)  any  $ \, K \in \HA_{\otimestilde} $  such that
$ \, K_0 := K \big/ h K \, $  is a topological Poisson Hopf algebra
isomorphic to  $ F[[\gerg]] $  for some finite dimensional Lie
bialgebra  $ \gerg $  (over  $ \Bbbk $);  then we write  $ \, H =
F_h[[\gerg]] \, $,  \, and say  $ K $  is a  {\sl quantisation\/}
of  $ F[[\gerg]] $.  We call \QFSHA{} the full subcategory of
$ \HA_{\otimestilde} $  whose objects are QFSHA, relative to all
possible  $ \gerg $  (see also  Remark 1.3{\it (a)\/}  below).
                                         \hfill\break
  \indent  (c) \, If  $ H_1 $,  $ H_2 $,  are two quantisations of
$ U(\gerg) $,  resp.~of  $ F[[\gerg]] $  (for some Lie bialgebra
$ \gerg $),  we say that  {\sl  $ H_1 $  is equivalent to  $ H_2 $},
and we write  $ \, H_1 \equiv H_2 \, $,  if there is an isomorphism
$ \, \varphi \, \colon H_1 \cong H_2 \, $  (in \QUEA, resp.~in \QFSHA)
%
%
%
%
such that  $ \, \varphi = \hbox{\it id} \mod h \, $.
\endproclaim

\vskip7pt

  {\bf Remarks 1.3:} \, {\it (a)} \, If  $ \, H \in \HA_{\otimeshat} $
is such that  $ \, H_0 := H \big/ h H \, $  as a Hopf algebra is
isomorphic to  $ U(\gerg) $  for some Lie algebra  $ \gerg $,  then
$ \, H_0 = U(\gerg) \, $  is also a  {\sl co-Poisson\/}  Hopf algebra
w.r.t.~the Poisson cobracket  $ \delta $  defined as follows: if  $ \,
x \in H_0 \, $  and  $ \, x' \in H \, $  gives  $ \, x = x' + h \, H
\, $,  \, then  $ \, \delta(x) := \big( h^{-1} \, \big( \Delta(x') -
\Delta^{\text{op}}(x') \big) \big) + h \, H \otimeshat H \, $;  \,
then (by [Dr], \S 3, Theorem 2) the restriction of  $ \delta $
makes  $ \gerg $  into a Lie bialgebra.  Similarly, if  $ \, K \in
\HA_{\otimestilde} $  is such that  $ \, K_0 := K \big/ h K \, $  is
a topological Poisson Hopf algebra isomorphic to  $ F[[\gerg]] $  for
some Lie algebra  $ \gerg $  then  $ \, K_0 = F[[\gerg]] \, $  is also
a topological  {\sl Poisson\/}  Hopf algebra w.r.t.~the Poisson bracket
$ \{\,\ ,\ \} $  defined as follows: if  $ \, x $,  $ y \in K_0 \, $
and  $ \, x' $,  $ y' \in K \, $  give  $ \, x = x' + h \, K $,  $ \,
y = y' + h \, K $,  \, then  $ \, \{x,y\} := \big( h^{-1} (x' \, y' -
y' \, x') \big) + h \, K \, $;  \, then  $ \gerg $  is a bialgebra
again, and  $ F[[\gerg]] $  is (the algebra of regular functions on)
a  {\sl Poisson\/}  formal group.  These natural co-Poisson and
Poisson structures are the ones considered in Definition 1.2 above.
                                                   \par
   In fact, specialisation gives a tensor functor from  \QUEA{}  to the
tensor category of universal enveloping algebras of Lie bialgebras and
a tensor functor from  \QFSHA{}  to the tensor category of (algebras
of regular functions on) formal Poisson groups.
                                                   \par
   {\it (b)} \, Clearly \QUEA, resp.~\QFSHA,  is a  {\sl
tensor\/}  subcategory of  $ \HA_{\otimeshat} $,  resp.~of
$ \HA_{\otimestilde} $.
                                                   \par
   {\it (c)} \, Let  $ H $  be a QFSHA.  Then  $ H $  is complete
w.r.t.~the weak topology, and  $ \, H_0 \cong F[[\gerg]] \, $  for
some finite dimensional Lie bialgebra  $ \gerg $,  and the weak
topology on  $ \, H_0 \cong F[[\gerg]] \, $  coincides with the
$ \hbox{\it Ker}\,(\epsilon_{\scriptscriptstyle H_0}) $--adic
topology.  It follows that the weak topology in  $ H $  coincides
with the  $ I_{\scriptscriptstyle H} $--adic  topology, so \QFSHA{}
is a subcategory of  $ \HA_{\otimestilde}^{w-I} $.  In particular,
if  $ \, H \in \QFSHA \, $  then  $ H \otimestilde H $  equals
the completion of  $ H \otimes_{\Bbbk[[h]]} H $  w.r.t.~the
$ I_{\scriptscriptstyle \! H \times H} $--adic  topology.

\vskip7pt

\proclaim {Definition 1.4}  Let  $ H $,  $ K $  be Hopf algebras (in
any category) over a ring  $ R $.  A pairing  $ \, \pi = \langle \,\ ,
\,\ \rangle \, \colon \, H \times K \loongrightarrow R \, $  is called
{\sl perfect\/}  if it is non-degenerate; it is called a  {\sl Hopf
pairing\/}  if for all  $ \, x $,  $ x_1 $,  $ x_2 \in H $,  $ y $,
$ y_1 $,  $ y_2 \in K $,  the elements  $ \, \big\langle \Delta(x),
y_1 \otimes y_2 \big\rangle := \sum_{(x)} \langle x_{(1)}, y_1 \rangle
\cdot \langle x_{(2)}, y_2 \rangle \, $  and  $ \, \big\langle x_1
\otimes x_2, \Delta(y) \big\rangle := \sum_{(y)} \langle x_1, y_{(1)}
\rangle \cdot \langle x_2, y_{(2)} \rangle \, $  are well defined
and we have
  $$  \displaylines{
   \big\langle x, y_1 \cdot y_2 \big\rangle = \big\langle \Delta(x),
y_1 \otimes y_2 \big\rangle \; ,  \qquad  \big\langle x_1 \cdot x_2,
y \big\rangle = \big\langle x_1 \otimes x_2, \Delta(y) \big\rangle
\cr
   \langle x, 1 \rangle = \epsilon(x) \; ,  \qquad  \langle 1,
y \rangle = \epsilon(y) \; , \qquad  \big\langle S(x), y \big\rangle
= \big\langle x, S(y) \big\rangle \; .  \cr }  $$
\endproclaim

\vskip7pt

  {\bf 1.5 Drinfeld's functors.} \,  Let  $ H $  be a Hopf
algebra (of any type) over  $ \Bbbk[[h]] $.  For each  $ \, n \in
\N $,  define  $ \; \Delta^n \colon H \longrightarrow H^{\otimes n}
\; $  by  $ \, \Delta^0 := \epsilon \, $,  $ \, \Delta^1 :=
{id}_{\scriptscriptstyle H} $,  \, and  $ \, \Delta^n := \big(
\Delta \otimes {id}_{\scriptscriptstyle H}^{\,\otimes (n-2)} \big)
\circ \Delta^{n-1} \, $  if  $ \, n \geq 2 $.  For any ordered subset
$ \, E = \{i_1, \dots, i_k\} \subseteq \{1, \dots, n\} \, $  with
$ \, i_1 < \dots < i_k \, $,  \, define the morphism  $ \;
j_{\scriptscriptstyle E} : H^{\otimes k} \longrightarrow
H^{\otimes n} \; $  by  $ \; j_{\scriptscriptstyle E} (a_1
\otimes \cdots \otimes a_k) := b_1 \otimes \cdots \otimes b_n
\; $  with  $ \, b_i := 1 \, $  if  $ \, i \notin \Sigma \, $
and  $ \, b_{i_m} := a_m \, $  for  $ \, 1 \leq m \leq k \, $;
then set
  $$  \Delta_E := j_{\scriptscriptstyle E} \circ \Delta^k \, ,  \quad
\Delta_\emptyset := \Delta^0 \, ,  \qquad  \text{and}  \qquad
\delta_E := \sum_{E' \subset E} {(-1)}^{n- \left| E' \right|}
\Delta_{E'} \, ,  \quad  \delta_\emptyset := \epsilon \, .  $$
By the inclusion-exclusion principle, the inverse formula  $ \;
\Delta_E = \sum_{\Psi \subseteq E} \delta_\Psi \, $  holds.  We
shall also use the notation  $ \, \delta_0 := \delta_\emptyset \, $,
$ \, \delta_n := \delta_{\{1, 2, \dots, n\}} \, $.  Then we define
  $$  H' := \big\{\, a \in H \,\big\vert\; \delta_n(a) \in
h^n H^{\otimes n} \; \forall\, n \in \N \,\big\}  \qquad
\big( \subseteq H \, \big ) \; .   $$
   \indent   Note that the useful formula  $ \; \delta_n =
{({id}_{\scriptscriptstyle H} - \epsilon)}^{\otimes n} \circ
\Delta^n \; $  holds, for all  $ \, n \in \N_+ \, $.  Then  $ H $
splits as  $ \, H = \Bbbk[[h]] \cdot 1_{\scriptscriptstyle H} \oplus
J_{\scriptscriptstyle H} \, $,  and  $ \, ({id} - \epsilon) \, $
projects  $ H $  onto  $ \, J_{\scriptscriptstyle H} := \text{\it
Ker}\, (\epsilon) \, $:  so  $ \, {({id} - \epsilon)}^{\otimes n}
\, $  projects  $ H^{\otimes n} $  onto  $ \, {J_{\scriptscriptstyle
H}}^{\! \otimes n} \, $;  therefore  $ \; \delta_n(a) = {({id}
- \epsilon)}^{\otimes n} \big( \Delta^n(a) \big) \in
{J_{\scriptscriptstyle H}}^{\! \otimes n} \; $  for
any  $ \, a \in H \, $.
                                        \par
   Now let  $ \, I_{\scriptscriptstyle H} := \epsilon^{-1}
\big( h \, \Bbbk[[h]] \big) \, $;  \, set  $ \, H^\times :=
\sum\limits_{n \geq 0} h^{-n} {I_{\scriptscriptstyle H}}^{\!n}
= \sum\limits_{n \geq 0} {\big( h^{-1} I_{\scriptscriptstyle H}
\big)}^n = \bigcup\limits_{n \geq 0} {\big( h^{-1}
I_{\scriptscriptstyle H} \big)}^n \, $  (the
$ \Bbbk[[h]] $--subalgebra  of  $ \lsF H $
generated by  $ \, h^{-1} I_{\scriptscriptstyle H}
\, $;  \, the second identity follows immediately from
$ \, {\big( h^{-1} I_{\scriptscriptstyle H} \big)}^n \subseteq
{\big( h^{-1} I_{\scriptscriptstyle H} \big)}^m \, $  for all
$ \, n < m \, $),  and define
  $$  H^\vee :=  h\text{--adic completion of the
$ \Bbbk[[h]] $--module }  H^\times  $$
({\sl Warning:\/}  $ H^\times $  naturally embeds into  $ \lsF H $,
whereas  $ H^\vee $  a priori  {\sl does not},  for the completion
procedure may "lead outside"  $ \lsF H \, $).  Note also that  $ \,
I_{\scriptscriptstyle H} = J_{\scriptscriptstyle H} + h \cdot H
\, $  (with  $ \, J_{\scriptscriptstyle H} \, $  as above), so
$ \, H^\times = \sum_{n \geq 0} h^{-n} {J_{\scriptscriptstyle
H}}^{\!n} \, $  and  $ \; H^\vee =  h\text{--adic completion of }
\sum\nolimits_{n \geq 0} h^{-n} {J_{\scriptscriptstyle H}}^{\!n} \, $.

\vskip7pt

   We are now ready to state the main result we are interested in:

\vskip7pt

\proclaim {Theorem 1.6} ("The quantum duality principle";
cf.~[Dr], \S 7)  Assume  $ \, \text{\it char\,}(\Bbbk) = 0 \, $.
                                        \hfill\break
   \indent   The assignments  $ \, H \mapsto H^\vee \, $  and
$ \, H \mapsto H' \, $  respectively define functors of tensor
categories  $ \, \QFSHA \longrightarrow \QUEA \, $  and  $ \, \QUEA
\longrightarrow \QFSHA \, $.  These functors are inverse to each
other.  Indeed,  for all  $ \, U_h(\gerg) \in \QUEA \, $  and all
$ \, F_h[[\gerg]] \in \QFSHA \, $  one has (cf.~\S 1.2)
  $$  {U_h(\gerg)}' \Big/ h \, {U_h(\gerg)}' = F[[\gerg^*]] \, ,
\qquad  {F_h[[\gerg]]}^\vee \Big/ h \, {F_h[[\gerg]]}^\vee =
U(\gerg^*)   \eqno \circledast  $$
that is,
   if  $ \, U_h(\gerg) $  is a quantisation of  $ \, U(\gerg) $  then
$ \, {U_h(\gerg)}' $  is a quantisation of  $ \, F[[\gerg^*]] $,  and
if  $ \, F_h[[\gerg]] $  is a quantisation of  $ \, F[[\gerg]] $  then
$ \, {F[[\gerg^*]]}^\vee $  is a quantisation of  $ \, U(\gerg^*) $.
%
%
%
%
Moreover,
the functors preserve equivalence, that is  $ \, H_1 \equiv H_2 \, $
implies  $ \, {H_1}^{\!\vee} \equiv {H_2}^{\!\vee} \, $  or
$ \, {H_1}' \equiv {H_2}' \, $.
\endproclaim

\vskip7pt

   Our analysis also move us to set the following (half-proved)

\vskip7pt

\proclaim {Conjecture 1.7} The quantum duality principle holds
as well for  $ \, \text{\it char\,}(\Bbbk) > 0 \, $.
\endproclaim

 \eject
%
%
%

\centerline {\bf \S \; 2 \  General properties of Drinfeld's functors }

\vskip10pt

   The rest of this paper will be devoted to prove Theorem 1.6.
In this section we establish some general properties of Drinfeld's
functors.  The first step is entirely standard.

\vskip7pt

\proclaim {Lemma 2.1}  The assignments  $ \, H \mapsto H^* = H^\star
\, $  and  $ \, H \mapsto H^\star \, $  define contravariant functors
of tensor categories  $ \; {(\ )}^\star \, \colon \,
{\Cal T}_{\otimeshat} \! \longrightarrow \! {\Cal P}_{\otimestilde}
\; $  and  $ \; {(\ )}^\star \, \colon \, {\Cal P}_{\otimestilde} \!
\longrightarrow \! {\Cal T}_{\otimeshat} \; $  which are inverse to
each other.  Their restriction gives antiequivalences of tensor
categories  $ \; \HA_{\otimeshat} \,{\buildrel {(\ )}^* \over
\loongrightarrow}\, \HA_{\otimestilde} \, $,  $ \; \HA_{\otimestilde}
\,{\buildrel {(\ )}^\star \over \loongrightarrow}\, \HA_{\otimeshat}
\, $,  \; and, if  $ \, \text{\it char}\, (\Bbbk) = 0 \, $,  $ \; \QUEA
\,{\buildrel {(\ )}^* \over \longrightarrow}\, \QFSHA \, $,  $ \;
\QFSHA \,{\buildrel {(\ )}^\star \over \longrightarrow}\, \QUEA \, $.
$ \square $
\endproclaim

\vskip7pt

   The following key fact shows that, in a sense, Drinfeld's
functors are  {\sl dual\/}  to each other:

\vskip7pt

\proclaim{Proposition 2.2}  Let  $ \, H \in \HA_{\otimeshat} \, $,
$ \, K \in \HA_{\otimestilde} \, $,  \, and let  $ \, \pi \! = \!
\langle \,\ , \,\ \rangle \, \colon H \times K \loongrightarrow
\Bbbk[[h]] \, $  be a Hopf pairing.  Then  $ \pi $  induces a bilinear
pairing $ \; \displaystyle{ \langle \,\ ,\,\ \rangle \, \colon H'
\! \times K^\vee \loongrightarrow \Bbbk[[h]] } \, $.
                                      \hfill\break
\indent   If in addition  $ \pi $ is perfect, and the induced
$ \Bbbk $--valued  pairing  $ \, \pi_0 \, \colon H_0 \times K_0
\loongrightarrow \Bbbk \, $  is still perfect, then  $ \; H' =
{\big( K^\times \big)}^\circ \! := \Big\{\, \eta \in \!
\lsF H \;\Big\vert\; \big\langle \eta, K^\times
\big\rangle \subseteq \Bbbk[[h]] \,\Big\} \; $  (w.r.t.~the
natural $ \Bbbk((h)) $--valued  pairing induced by scalar extension).
In particular, if  $ \, H = K^\star \, $  and  $ \, K = H^* \, $
the evaluation pairing yields  $ \Bbbk[[h]] $--module  isomorphisms
$ \; H' \cong {\big( K^\vee \big)}^* \, $  and  $ \; K^\vee \cong
{\big( H' \big)}^\star \, $.
\endproclaim

\demo{Proof}  First note that, for all  $ \, x $,  $ x_1 $,  $ x_2
\in H $,  $ y $,  $ y_1 $,  $ y_2 \in K $,  the elements  $ \,
\, \big\langle \Delta(x), y_1 \otimes y_2 \big\rangle := \sum_{(x)}
\langle x_{(1)}, y_1 \rangle \cdot \langle x_{(2)}, y_2 \rangle \, $
and  $ \, \big\langle x_1 \otimes x_2, \Delta(y) \big\rangle :=
\sum_{(y)} \langle x_1, y_{(1)} \rangle \cdot \langle x_2, y_{(2)}
\rangle \, $  (cf.~Definition 1.4) are well defined: in fact  $ K $
acts   --- via  $ \pi $  ---   as a Hopf subalgebra of  $ H^* $,
hence  $ \, K \otimestilde K \, $  acts   --- via  $ \, \pi \otimes
\pi \, $  ---   as a Hopf subalgebra of  $ \, H^* \otimestilde H^* =
{\big( K \otimeshat K \big)}^* \, $,  \, due to Lemma 2.1.  Therefore
it is perfectly meaningful to require $ \pi $  to be a  {\sl Hopf\/}
pairing.
                                                 \par
   Now, scalar extension gives a Hopf pairing  $ \, \langle \,\ ,\,
\ \rangle \, \colon \lsF H \times \lsF K \longrightarrow \Bbbk((h)) \, $
which restricts to a similar pairing  $ \, \langle \,\ ,\,\ \rangle
\, \colon \, H' \times K^\times \longrightarrow \Bbbk((h)) \, $:  we have
to prove that the latter takes values in  $ \, \Bbbk[[h]] \, $,  that is
$ \, \Big\langle H', K^\times \Big\rangle \subseteq \Bbbk[[h]] \, $,  for
then it will extend by continuity to a pairing  $ \, \langle \,\ , \,
\ \rangle \, \colon \, H' \times K^\vee \longrightarrow \Bbbk[[h]] \, $;
in addition, this will also imply  $ \, H' \subseteq {\big( K^\times
\big)}^\circ \, $.
                                                 \par
   Take  $ \, c_1 $,  $ \dots $,  $ c_n \in I_{\scriptscriptstyle K}
\, $;  then  $ \, \langle 1, c_i \rangle = \epsilon(c_i) \in h \,
\Bbbk[[h]] \, $.  Now, given  $ \, y \in H' $,  \, consider
  $$  \left\langle y \, , \, \prod_{i=1}^n c_i \right\rangle =
\left\langle \Delta^n(y) \, , \, \otimes_{i=1}^n c_i \right\rangle
= \left\langle \sum_{\Psi \subseteq \{1,\dots,n\}} \hskip-11pt
\delta_\Psi(y) \, , \, \otimes_{i=1}^n c_i \right\rangle =
\sum_{\Psi \subseteq \{1,\dots,n\}} \hskip-9pt \big\langle
\delta_\Psi(y) \, , \, \otimes_{i=1}^n c_i \big\rangle  $$
(using formulas in \S 1.5) and look at the generic summand in
the last expression above.  Let  $ \, \vert\Psi\vert = t \, $
($ \, t \leq n \, $):  then  $ \, \big\langle \delta_\Psi(y) ,
\otimes_{i=1}^n c_i \big\rangle = \big\langle \delta_t(y) ,
\otimes_{i \in \Psi} c_i \big\rangle \cdot \prod_{j \not\in \Psi}
\langle 1, c_j \rangle \, $,  by definition of  $ \, \delta_\Psi \, $.
Thanks to the previous analysis, we have  $ \, \prod_{j \not\in \Psi}
\langle 1, c_j \rangle \in h^{n-t} \Bbbk[[h]] \, $,  and  $ \, \big\langle
\delta_t(y), \otimes_{i \in \Psi} c_i \big\rangle \in h^t \Bbbk[[h]] \, $
because  $ \, y \in H' \, $;  thus we get  $ \, \big\langle \delta_t(y),
\otimes_{i \in \Psi} c_i \big\rangle \cdot \prod_{j \not\in \Psi}
\langle 1, c_j \rangle \in h^n \Bbbk[[h]] \, $,  whence  $ \, \big\langle
y, \prod_{i=1}^n c_i \big\rangle \in h^n \Bbbk[[h]] \, $.  The outcome is
that  $ \; \displaystyle{ \langle \, y, \psi \rangle \in h^n \Bbbk[[h]] }
\; $  for all  $ \, y \in H' $,  $ \psi \in {I_{\scriptscriptstyle
K}}^{\hskip-2pt n} $,  \, and therefore  $ \, \big\langle H', \,
h^{-n} {I_{\scriptscriptstyle K}}^{\hskip-2pt n} \big\rangle \subseteq
\Bbbk[[h]] \, $  for all  $ n \in \N $,  \, whence  $ \, \big\langle H',
K^\times \big\rangle \subseteq \Bbbk[[h]] \, $,  \, q.e.d.
                                             \par
   We are now left with proving  $ \, {\big( K^\times \big)}^\circ
\subseteq H' \, $:  we do it by reverting the previous argument.
                                             \par
   Let  $ \, \eta \in {\big( K^\times \big)}^\circ \, $:  then
$ \, \big\langle \eta \, , \, h^{-s} {I_{\scriptscriptstyle
K}}^{\hskip-2pt s} \big\rangle \in \Bbbk[[h]] \, $  hence  $ \,
\big\langle \eta \, , {I_{\scriptscriptstyle K}}^{\hskip-2pt s}
\big\rangle \in h^s \Bbbk[[h]] \, $,  for all  $ \, s \in \N \, $.
In particular, for  $ \, s = 0 \, $  this gives  $ \, \big\langle
\eta \, , K \big\rangle \in \Bbbk[[h]] \, $,  whence   --- thanks to
non-degeneracy of  $ \pi_0 $  ---   we get $ \, \eta \in H \, $.
Let now  $ \, n \in \N \, $  and  $ \, i_1 $,  $ \dots $,
$ \, i_n \in I_{\scriptscriptstyle K} \, $; then
  $$  \displaylines{
   {} \hskip7pt   \Big\langle \delta_n(\eta) \, , \, \otimes_{k=1}^n
i_k \Big\rangle \, = \, \left\langle \sum\nolimits_{\Psi \subseteq
\{1,\dots,n\}} {(-1)}^{n - \vert \Psi \vert} \, \Delta_\Psi(\eta)
\, , \, \otimes_{k=1}^n i_k \right\rangle \, =   \hfill {}  \cr
   = \, \sum\nolimits_{\Psi \subseteq \{1,\dots,n\}} {(-1)}^{n -
\vert \Psi \vert} \cdot \left\langle \eta \, ,
\, \prod\nolimits_{k \in \Psi} i_k \right\rangle \cdot
\prod\nolimits_{k \not\in \Psi} \langle 1, i_k \rangle \, \in
\hskip55pt {}  \cr
%
%
   {} \hfill   \in \, \sum\nolimits_{\Psi \subseteq \{1,\dots,n\}}
\big\langle \eta \, , \, {I_{\scriptscriptstyle K}}^{\hskip-2pt
\vert \Psi \vert} \big\rangle \cdot h^{n - \vert \Psi \vert} \,
\Bbbk[[h]] \; \subseteq \; \sum\nolimits_{s=0}^n \, h^s \cdot h^{n-s}
\, \Bbbk[[h]] = h^n \, \Bbbk[[h]] \, ,  \cr }  $$
therefore  $ \, \Big\langle \delta_n(\eta) \, , \,
{I_{\scriptscriptstyle K}}^{\hskip-3pt \otimes n} \Big\rangle
\subseteq h^n \Bbbk[[h]] \, $.  In addition,  $ H $  splits as
$ \, K = \Bbbk[[h]] \cdot 1_{\scriptscriptstyle K} \oplus
J_{\scriptscriptstyle K} \, $,  so  $ K^{\otimes n} $  splits
into the direct sum of  $ \, {J_{\scriptscriptstyle K}}^{\hskip-3pt
\otimes n} \, $  and of other direct summands which are again tensor
products but in which at least one tensor factor is  $ \, \Bbbk[[h]] \cdot
1_{\scriptscriptstyle K} \, $.  Since  $ \; J_{\scriptscriptstyle H}
:= \hbox{\it Ker}\, (\epsilon_{\scriptscriptstyle H}) = {\big\langle
\Bbbk[[h]] \cdot 1_{\scriptscriptstyle K} \big\rangle}^\perp = \big\{\,
y \in H \,\big\vert\, \langle \, y, 1_{\scriptscriptstyle K} \rangle
= 0 \,\big\} \; $  (the subspace of  $ H $  orthogonal to  $ \,
\big\langle \Bbbk[[h]] \cdot 1_{\scriptscriptstyle K} \big\rangle \, $),
we have  $ \, \Big\langle {J_{\scriptscriptstyle H}}^{\hskip-3pt
\otimes n}, \, K^{\otimes n} \Big\rangle = \Big\langle
{J_{\scriptscriptstyle H}}^{\hskip-3pt \otimes n}, \,
{J_{\scriptscriptstyle K}}^{\hskip-3pt \otimes n}
\Big\rangle \, $.  Since  $ \, \delta_n(\eta) \in
{J_{\scriptscriptstyle H}}^{\hskip-3pt \otimes n} \, $
(cf.~\S 1.5), this analysis yields  $ \, \Big\langle \delta_n(\eta)
\, , K^{\otimes n} \Big\rangle \subseteq h^n \, \Bbbk[[h]] \, $,  whence
--- due to the non-degeneracy of the specialised pairing ---   we get
$ \, \delta_n(\eta) \in h^n H^{\otimes n} \, $.  Therefore  $ \, \eta
\in H' $;  hence we get  $ \, {\big( K^\times \big)}^\circ \subseteq
H' $,  q.e.d.
                                             \par
   For the last part of the statement, since  $ K^\vee $  is the
$ h $--adic  completion of  $ K^\times $  one has  $ \, {\big( K^\vee
\big)}^* = {\big( K^\times \big)}^* \, $,  \, so now we show that the
latter is equal to  $ \, {\big( K^\times \big)}^\circ = H' \, $.  On
the one hand, it is clear that  $ \, H' = {\big( K^\times \big)}^\circ
\subseteq {\big( K^\times \big)}^* $.  On the other hand, pick  $ \,
f \in {\big( K^\times \big)}^* \, $:  \, then  $ f $  is uniquely
determined by  $ \, f{\Big\vert}_K \, $,  \, and by construction
$ \, f{\Big\vert}_K \in K^* \, $  and  $ \, f{\Big\vert}_K \big(
{I_{\scriptscriptstyle K}}^{\hskip-3pt n} \big) \subseteq h^n \Bbbk[[h]]
\, $  because  $ \, f \big( h^{-n} {I_{\scriptscriptstyle K}}^{\hskip
-3pt n} \big) \subseteq f \big( K^\times \big) \subseteq \Bbbk[[h]] \, $.
Therefore  $ \, f{\Big\vert}_K \in K^\star = {\big( H^* )}^\star =
H \, $  (by Lemma 2.1), thus  $ \, f{\Big\vert}_K \! \in H \, $  and
$ \, f{\Big\vert}_K \big( K^\times \big) \subseteq \Bbbk[[h]] \, $
yields  $ \, f{\Big\vert}_K \! \in {\big( K^\times \big)}^\circ
\! = H' \, $,  \, whence  $ \, f \in H' $.   $ \square $
\enddemo

\vskip7pt

\proclaim{Lemma 2.3} Let  $ \, H_1, H_2 \in \HA_{\otimestilde}^{w-I}
\, $.  Then  $ \, {\big( H_1 \otimestilde H_2 \big)}^{\!\vee} =
{H_1}^{\!\vee} \otimeshat {H_2}^{\!\vee} \, $.  In particular
this holds true for any  $ \, H_1, H_2 \in \QFSHA \, $.
\endproclaim

\demo{Proof} Clearly  $ \, I_{\scriptscriptstyle \! H_1 \!\otimestilde\!
H_2} = I_{\scriptscriptstyle \! H_1} \otimestilde H_2 + H_1 \otimestilde
I_{\scriptscriptstyle \! H_2} \, $,  \, and the assumption on topologies
implies that  $ H_1 \otimestilde H_2 $  is the  $ I_{\scriptscriptstyle
\! H_1 \!\otimestilde\! H_2} $--adic  completion of  $ H_1
\otimes_{\Bbbk[[h]]} H_2 $.  Then, for each  $ \, \eta_{\otimestilde}
\in H \otimestilde H \, $  we can find an expression  $ \,
\eta_{\otimestilde} = \sum_{m \in \N} \eta_{\otimestilde}^{(m)}
\, $  such that  $ \, \eta_{\otimestilde}^{(m)} \in
{(I_{\scriptscriptstyle \! H \!\otimes\! H})}^m $  for
all  $ m \, $;  as  $ {(I_{\scriptscriptstyle \! H
\!\otimestilde\! H})}^m $  is the completion  $ \,
\sum_{r+s=m} I^r \otimestilde I^s $  of  $ \, \sum_{r+s=m}
I^r \otimes I^s $,  we can in fact write  $ \, \eta_{\otimestilde}
= \sum_{m \in \N} \sum_{r+s=m} \eta_m^{(r)} \otimes \eta_m^{(s)} \, $
for some  $ \, \eta_m^{(r)} \in I^r $,  $ \, \eta_m^{(s)} \in I^s $
(for all  $ m $,  $ r $,  $ s $),  with  $ \, \sum_{r+s=m} \eta_m^{(r)}
\otimes \eta_m^{(s)} = 0 \, $  for all  $ \, m < n \, $  if  $ \,
\eta_{\otimestilde} \in {I_{\scriptscriptstyle \! H \!\otimestilde\!
H}}^{\hskip-3pt n} \, $.  Thus for any  $ \, n \! \in \! \N \, $
and  $ \, \eta_{\otimestilde} \! \in \! {I_{\scriptscriptstyle \!
H \!\otimestilde\! H}}^{\hskip-3pt n} \, $
  $$  \displaylines{
   {} \quad   h^{-n} \eta_{\otimestilde} = h^{-n} \sum_{m \geq n}
\sum_{r+s=m} \eta_m^{(r)} \otimes \eta_s^{(s)} \in h^{-n} \,
\sum_{m \geq n} \sum_{r+s=m} I^r \otimes I^s =   \hfill   \cr
   \hfill   = \sum_{m \geq n} \sum_{r+s=m} \, h^{m-n} \, h^{-r}
I^r \otimes h^{-s} I^s \subseteq \sum_{\ell \in \N} h^\ell \,
H^\times \otimes H^\times \, ,   \quad  \cr }  $$
from which one argues that the natural morphism  $ \, H^\times \!
\loongrightarrow H^\vee \, $  induces a similar map  $ \, {\big(
H \otimestilde H \big)}^\vee \! \loongrightarrow H^\vee \otimeshat
H^\vee \, $.  Conversely,  $ \, \sum_{r+s=m} {I_{\scriptscriptstyle
\! H_1}}^{\hskip-2pt r} \otimes {I_{\scriptscriptstyle \!
H_2}}^{\hskip-2pt s} \subseteq {I_{\scriptscriptstyle \! H_1
\!\otimestilde\! H_2}}^{\hskip-3pt m} \, $  (for all  $ m $)
implies  $ \, {H_1}^{\!\times} \!\otimes {H_2}^{\!\times}
\subseteq {\big( H_1 \!\otimestilde\! H_2 \big)}^\times $,
\, whence one gets by completion a continuous morphism  $ \,
{H_1}^{\!\vee} \!\otimeshat {H_2}^{\!\vee} \! \longrightarrow
{\big( H_1 \!\otimestilde\! H_2 \big)}^{\!\vee} $,  \, inverse
of the previous one.  This gives the equality in the claim.
                                        \par
   Finally, by Remark 1.3{\it(c)}  any  $ \, H_1, H_2 \in \QFSHA \, $
fulfills the hypotheses.   $ \square $
\enddemo

\vskip7pt

\proclaim{Proposition 2.4}
                                  \hfill\break
   \indent   (a) \, Let  $ \, H \in \HA_{\otimestilde} \, $.
Then  $ \, H^\vee \, $  is a unital (topological)  $ \Bbbk[[h]] $--algebra
in  $ {\Cal T}_{\otimeshat} \, $.
                                  \hfill\break
   \indent   (b) \, Let  $ \, H \in \HA_{\otimestilde}^{w-I} \, $.
Then  $ \, H^\vee \in \HA_{\otimeshat} \, $,  \, and the  $ \Bbbk $--Hopf
algebra  $ \, {H^\vee}_{\hskip-3pt 0} \, $  is cocommutative and
connected; if  $ \, \text{\it char\,}(\Bbbk) = 0 \, $,  \, it is a
universal enveloping algebra, and  $ \, H^\vee \in \QUEA \, $.
\endproclaim

\demo{Proof}  {\it (a)} \, We must prove that  $ H^\vee $  is
topologically free: by the criterion in [KT], \S 4.1, this is
equivalent to  $ H^\vee $  being a torsionless, separated and
complete topological  $ \Bbbk[[h]] $--module.  Now,  $ H $  is
torsionless, so the same is true for  $ \lsF H $  hence
for  $ H^\times $  too; as  $ H^\vee $  is the  $ h $--adic
completion of  $ H^\times $,  it is torsionless as well, and
by definition it is complete and separated.  Furthermore, by
construction  $ H^\vee $  is a (topological)  $ \Bbbk[[h]] $--algebra,
unital since  $ \, 1_{\scriptscriptstyle H} \in H^\times $.
                                        \par
   {\it (b)} \, Let  $ \, I := I_H \, $  (cf.~\S 1.5).  The definition
yields  $ \, S_H(I) = I \, $,  whence  $ \, S_H \big( h^{-n} I^n \big)
= h^{-n} I^n \, $  for all  $ \, n \in \N $,  so  $ \, S_H \big(
H^\times \big) = H^\times $,  \, so one can define  $ \, S_{H^\vee}
\, $  by continuous extension.  As for  $ \Delta $,  the assumption on
topologies implies that  $ \, {\big( H \otimestilde H \big)}^\vee =
H^\vee \otimeshat H^\vee $,  \, by Lemma 2.3.  Moreover, definitions
yield  $ \, \Delta_H \big( I^n \big) \subseteq \sum_{r+s=n} I^r
\otimestilde I^s = {I_{\scriptscriptstyle \! H \!\otimestilde\!
H}}^{\hskip-3pt n} $  (for all  $ n \, $),  hence
  $$  \Delta_H \big( h^{-n} I^n \big) \subseteq h^{-n} \sum_{r+s=n}
I^r \otimestilde I^s = h^{-n} {I_{\scriptscriptstyle \! H
\!\otimestilde\! H}}^{\hskip-3pt n} \subseteq {\big( H \otimestilde
H \big)}^{\!\times}  $$
so that  $ \, \Delta_H \big( H^\times \big) \subseteq {\big(
H \otimestilde H \big)}^{\!\times} $,  \, thus one gets  $ \,
\Delta_{H^\vee} \, $  by continuity.  Finally, by construction
$ \epsilon_H $  extends to a counit for  $ H^\vee $.  It is clear
that all axioms of a Hopf algebra in  $ {\Cal T}_{\otimeshat} $
are then fulfilled, therefore  $ \, H^\vee \in \HA_{\otimeshat}
\, $.  Now, since  $ \, H^\times = \sum_{n \geq 0} {\big( h^{-1}
J_{\scriptscriptstyle H} \big)}^n \, $,  \, the unital
topological algebra  $ H^\vee $  is generated by  $ \,
{J_{\scriptscriptstyle H}}^{\!\!\times} \! := \! h^{-1}
J_{\scriptscriptstyle H} \, $.  Consider  $ \, j^\vee \!
\in \! {J_{\scriptscriptstyle H}}^{\!\!\times} \, $,  \, and
$ \, j := h \, j^\vee \! \in \! J_{\scriptscriptstyle H} \, $;
\, then  $ \; \Delta = \delta_2 + id \otimes 1 + 1 \otimes id -
\epsilon \cdot 1 \otimes 1 \; $  and  $ \, \hbox{\it Im}\,(\delta_2)
\subseteq J_{\scriptscriptstyle H} \otimes J_{\scriptscriptstyle H}
\, $  (cf.~\S 1.5) give
  $$  \Delta(j) = \delta_2(j) + j \otimes 1 + 1 \otimes j -
\epsilon(j) \cdot 1 \otimes 1 \in j \otimes 1 + 1 \otimes j +
J_{\scriptscriptstyle H} \otimestilde J_{\scriptscriptstyle H}
\, .  $$
   Therefore
  $$  \Delta \big( j^\vee \big) = \delta_2 \big( j^\vee \big) +
j^\vee \otimes 1 + 1 \otimes j^\vee - \epsilon \big(j^\vee\big)
\, 1 \otimes 1 = j^\vee \otimes 1 + 1 \otimes j^\vee + \delta_2
\big( j^\vee \big)  $$
which maps (through completion) into
  $$  j^\vee \otimes 1 + 1 \otimes j^\vee +
h^{-1} J_{\scriptscriptstyle H} \otimestilde
J_{\scriptscriptstyle H} = j^\vee \otimes 1 +
1 \otimes j^\vee + h^{+1} {J_{\scriptscriptstyle H}}^{\!\!\vee}
\otimeshat {J_{\scriptscriptstyle H}}^{\!\!\vee} \, ,  $$
whence we conclude that
  $$  {} \qquad  \hfill   \Delta_{H^\vee} \big( j^\vee \big)
\equiv j^\vee \otimes 1 + 1 \otimes j^\vee \mod h \, H^\vee
\otimeshat H^\vee   \hskip45pt   \forall \;\; j^\vee \in
{J_{\scriptscriptstyle H}}^{\!\!\vee} \, .  \qquad  $$
Thus  $ \, {J_{\scriptscriptstyle H}}^{\!\!\vee} \mod h \, H^\vee
\, $ is contained in  $ P \big( H^\vee \big) $,  the set of  {\sl
primitive elements\/}  of  $ {H^\vee}_{\!\!0} \, $;  since  $ \,
{J_{\scriptscriptstyle H}}^{\!\!\vee} \mod h \, H^\vee \, $  generates
$ \, {H^\vee}_{\!\!0} \, $   --- as  $ \, {J_{\scriptscriptstyle
H}}^{\!\!\vee} \, $  generates  $ \, H^\vee $  ---   this proves
{\it a fortiori}  that  $ P \big( H^\vee \big) $  generates  $ \,
{H^\vee}_{\!\!0} \, $,  and also shows that  $ \, {H^\vee}_{\!\!0}
\, $  is cocommutative.  In addition, we can also apply Lemma 5.5.1
in [M] to the Hopf algebra  $ \, {H^\vee}_{\!\!0} \, $,  with  $ \,
A_0 = \Bbbk \cdot 1 \, $  and  $ \, A_1 = {J_{\scriptscriptstyle H}}^{\!
\!\vee} \mod h \, H^\vee $,  \, to argue that  $ \, {H^\vee}_{\!\!0}
\, $  is connected, q.e.d.
                                        \par
   If  $ \, \text{\it char\,}(\Bbbk) = 0 \, $  by Kostant's Theorem
(cf.~for instance [A], Theorem 2.4.3) we have  $ \, {H^\vee}_{\!\!0}
= U(\gerg) \, $  for the Lie (bi)algebra  $ \, \gerg = P \big(
{H^\vee}_{\!\!0} \big) \, $.  We conclude that  $ \, H \in
\QUEA \, $.   $ \square $
\enddemo

\vskip7pt

\proclaim{Lemma 2.5} ([KT], Lemma 3.2) Let  $ H $  be a Hopf
$ \Bbbk[[h]] $--algebra,  let  $ \, a $,  $ b \in H $,  and let
$ \Phi $  be a finite subset of  $ \, \N \, $.  Then  $ \;
\displaystyle{ \delta_\Phi(ab) = {\textstyle \sum}_{ \Lambda \cup
Y = \Phi} \delta_\Lambda(a) \, \delta_Y(b) } \, $.  In addition,
if  $ \, \Phi \not= \emptyset \, $  then  $ \; \displaystyle{
\delta_\Phi(ab - ba) = {\textstyle \sum}_{\Sb  \Lambda \cup Y = \Phi
\\   \Lambda \cap Y \not= \emptyset  \endSb}  \big( \delta_\Lambda(a)
\, \delta_Y(b) - \delta_Y(b) \, \delta_\Lambda(a) \big) } \, $.
$ \square $
\endproclaim

\vskip7pt

\proclaim{Proposition 2.6} Let  $ \, H \in \HA_{\otimeshat} \, $.
Then  $ \, H' \in \HA_{\otimestilde} \, $,  \, and the  $ \Bbbk $--Hopf
algebra  $ \, {H'}_{\hskip-1pt 0} \, $  is commutative.
\endproclaim

\demo{Proof} First,  $ H' $  is a  $ \Bbbk[[h]] $--submodule  of  $ H $,
for the maps  $ \delta_n $  ($ n \in \N $)  are  $ \Bbbk[[h]] $--linear;
to see it lies in  $ {\Cal P}_{\otimestilde} $,  we resort to a duality
argument.  Let  $ \, K := H^* \in \HA_{\otimestilde} \, $,  \, so  $ \,
H = K^\star \, $  (cf.~Lemma 2.1), and let  $ \, \pi \, \colon H \!
\times \! K \! \loongrightarrow \Bbbk[[h]] \, $  be the natural Hopf
pairing given by evaluation.  Then Proposition 2.2 gives  $ \, H'
\cong {\big( K^\vee \big)}^* \! \in \! {\Cal P}_{\otimestilde} \, $,
\, thus since  $ K^\vee $  is a unital algebra we have that  $ H' $
is a counital coalgebra in  $ {\Cal P}_{\otimestilde} $,  with  $ \,
\Delta_{\scriptscriptstyle H'} = \Delta_{\scriptscriptstyle H}
\!\big\vert_{\scriptscriptstyle H'} \, $  and  $ \,
\epsilon_{\scriptscriptstyle H'} = \epsilon_{\scriptscriptstyle H}
\!\big\vert_{\scriptscriptstyle H'} \, $.  In addition, by Lemma 2.5
one easily sees that  $ H' $  is a  $ \Bbbk[[h]] $--subalgebra  of  $ H $,
and by construction it is unital for  $ \, 1_{\scriptscriptstyle H}
\in H' \, $.  The outcome is  $ \, H' \in \HA_{\otimestilde} \, $,
\, q.e.d.
                                        \par
   Finally, the very definitions give  $ \, x = \delta_1(x) +
\epsilon(x) \, $  for all  $ \, x \in H \, $.  If  $ \, x \in H'
\, $  we have  $ \, \delta_1(x) \in h \, H \, $,  hence there exists
$ \, x_1 \in H \, $  such that  $ \, \delta_1(x) = h x_1 \, $.  Now
for  $ \, a, b \in H' $,  \, write  $ \; a = h \, a_1 + \epsilon(a)
\, $,  $ \; b = h \, b_1 + \epsilon(b) \, $,  \, hence  $ \; ab - ba
= h \, c \; $  with  $ \; c = h \, (a_1 b_1 - b_1 a_1) \, $;  \,
we show that  $ \, c \in H' $.  For this we must check that
$ \delta_\Phi(c) $  is divisible by  $ h^{\vert \Phi \vert} $
for any finite subset  $ \Phi $  of  $ \N_+ \, $:  as multiplication
by  $ h $  is injective (for  $ H $  is topologically free), it is
enough to show that  $ \delta_\Phi(a b - b a) $  is divisible by
$ h^{\vert \Phi \vert + 1} $.  Let  $ \Lambda $  and  $ Y $  be
subsets of  $ \Phi $  such that  $ \, \Lambda \cup Y = \Phi \, $
and  $\, \Lambda \cap Y \not= \emptyset \, $:  then  $ \, \vert
\Lambda \vert + \vert Y \vert \geq \vert \Phi \vert + 1 \, $.
Now,  $ \delta_\Lambda(a) $  is divisible by  $ h^{\vert \Lambda
\vert} $  and  $ \delta_Y(b) $  is divisible by  $ h^{\vert Y \vert}
\, $:  \, from this and the second part of Lemma 2.5 it follows that
$ \delta_\Phi(ab - ba) $  is divisible by  $ h^{\vert \Phi \vert
+ 1} $.  The outcome is  $ \; ab \equiv ba \mod h \, H' \, $,  \,
so  $ {H'}_{\hskip-1pt 0} $  is commutative.   $ \square $
\enddemo

\vskip7pt

\proclaim{Lemma 2.7} Let  $ \, H_1 $,  $ H_2 \in \QUEA \, $.
Then  $ \, {\big( H_1 \otimeshat H_2 \big)}' = {H_1}'
\otimestilde {H_2}' \, $.
\endproclaim

\demo{Proof}  Proceeding as in the proof of Proposition 2.6, let
$ \, K_i := {H_i}^* \in \QFSHA \, $  ($ i = 1, 2 $);  then  $ \,
K_1 \otimestilde K_2 = {H_1}^* \otimestilde {H_2}^* = {\big( H_1
\otimeshat H_2 \big)}^* \, $  (by Lemma 2.1),  and  $ \, {H_i}' =
{\big( {K_i}^{\!\vee} \big)}^* \, $  ($ i = 1, 2 $),  and similarly
$ \, {\big( H_1 \otimeshat H_2 \big)}' \! = \! \Big( \! {\big( K_1
\otimestilde K_2 \big)}^{\!\vee} \Big)^{\!*} $.  Then applying Lemma
2.3 we get  $ \, {\big( H_1 \otimeshat H_2 \big)}' \! = \Big( \! {\big(
K_1 \otimestilde K_2 \big)}^{\!\vee} \Big)^{\!*} = {\big( {K_1}^{\!\vee}
\otimeshat {K_2}^{\!\vee} \big)}^* = {\big( {K_1}^{\!\vee} \big)}^*
\otimestilde {\big( {K_2}^{\!\vee} \big)}^* = {H_1}' \otimestilde
{H_2}' $,  \, q.e.d.   $ \square $
\enddemo

\vskip7pt

\proclaim{Lemma 2.8}  The assignment  $ \, H \mapsto H^\vee \, $,
\, resp.~$ H \mapsto H' \, $,  \, gives a well-defined functor
$ \, \HA_{\otimestilde}^{w-I} \! \loongrightarrow \HA_{\otimeshat}
\, $,  \, resp.~$ \, \HA_{\otimeshat} \loongrightarrow
\HA_{\otimestilde} \, $.
\endproclaim

\demo{Proof}  In order to define the functors, we only have to
set them on morphisms.  Let  $ \, H, K \in \HA_{\otimestilde}^{w-I}
\, $  and  $ \, \phi \in \hbox{\it Mor}_{\HA_{\otimestilde}^{w-I}}
(H,K) \, $;  \, by scalar extension it gives a morphism  $ \, \lsF
H \longrightarrow \! \lsF K \, $  of  $ \Bbbk((h)) $--Hopf algebras,
which maps  $ \, h^{-1} I_{\scriptscriptstyle H} \, $  into
$ \, h^{-1} I_{\scriptscriptstyle K} \, $,  hence  $ H^\times $
into  $ K^\times \, $:  \, extending it by continuity we get
the morphism  $ \, \phi^\vee \in \hbox{\it Mor}_{\scriptscriptstyle
\HA_{\otimeshat}} \big( H^\vee, K^\vee \big) \, $  we were looking
for.  Similarly, let  $ \, H, K \in \HA_{\otimeshat} \, $  and  $ \,
\varphi \in \hbox{\it Mor}_{\scriptscriptstyle \HA_{\otimeshat}}(H,K)
\, $:  \, then  $ \, \delta_n \circ \varphi = \varphi^{\otimes n}
\circ \delta_n \, $  (for all  $ n \in \N $),  so  $ \, \varphi
\big( H' \big) \subseteq K' \, $:  \, thus as  $ \, \varphi' \in
\hbox{\it Mor}_{\scriptscriptstyle \HA_{\otimestilde}} \big( H',
K' \big) \, $  we simply take  $ \, \varphi{\big\vert}_{H'} \, $.
$ \square $
\enddemo

\vskip1,5truecm

\centerline {\bf \S \; 3 \  Drinfeld's functors on quantum groups }

\vskip10pt

   We focus now on the effect of Drinfeld's functors on
quantum groups.  The first result is an explicit description
of  $ {F_h[[\gerg]]}^\vee $  when  $ F_h[[\gerg]] $  is a QFSHA.

\vskip7pt

   {\bf 3.1 An explicit description of  $ {F_h[[\gerg]]}^\vee $.} \,
Let  $ \, F_h[[\gerg]] \in \QFSHA \, $,  \, and set for simplicity
$ \, F_h := F_h[[\gerg]] \, $,  $ \, F_0 := F_h \big/ h \, F_h =
F[[\gerg]] \, $,  $ \, {F_h}^{\!\!\vee} := {F_h[[\gerg]]}^\vee $,
\, and  $ \, {F_0}^{\!\!\vee} := {F_h}^{\!\!\vee} \big/ h \,
{F_h}^{\!\!\vee} $.  Then  $ \, F_0 \cong F[[\gerg]] = \Bbbk[[
\bar{x}_1, \dots, \bar{x}_n ]] \, $  (for some  $ \, n \in
\N \, $)  as topological  $ \Bbbk $--algebras.  Letting  $ \, \pi
\colon \, F_h \relbar\joinrel\twoheadrightarrow F_0 \, $  be
the natural projection, if we pick an  $ \, x_j \in
       \pi^{-1}(\bar{x}_j) \, $  for any  $ j $,\break
         then  $ F_h $  is generated by  $ \, \{x_1,\dots,x_n\} \, $
as a topological  $ \Bbbk[[h]] $--algebra,  that is to say  $ \; F_h =
F_h[[\gerg]] = \Bbbk[[x_1,\dots,x_n,h]] \, $.  In this description we
have  $ \, I := I_{\scriptscriptstyle F_h} = \big(x_1,\dots,x_n,h\big)
\, $,  \, and  $ I^\ell $  identifies with the space of all formal
series whose degree (that is, the degree of the lowest degree
monomials occurring in the series with non-zero coefficient) is
at least  $ \ell $,  that is
  $$  I^\ell = \Big\{\, f = \sum\nolimits_{\bold{d} \in \N^{n+1}}
c_{\bold{d}} \cdot h^{d_0} x_1^{d_1} x_2^{d_2} \cdots x_n^{d_n}
\;\Big|\;\, c_{\bold{d}} \in \Bbbk, \, \bold{d} \in \N^{n+1}, \;
|\bold{d}| \geq \ell \,\Big\}  $$
for all  $ \, \ell \in \N $  (hereafter, we set  $ \, |\bold{d}| :=
\sum_{s=0}^n d_s \, $  for any  $ \, \bold{d} = (d_0,d_1,\dots,d_n)
\in \N^{n+1} \, $).  Then  $ \lsF F_h \cong \Bbbk[[x_1, \dots, x_n]]((h))
\, $,  \, and
  $$  h^{-\ell} I^\ell = \Big\{\, f = \sum\nolimits_{\bold{d} \in
\N^{n+1}} \! P_n(\underline{\check{x}}) \cdot h^n \;\Big|\;
P_n(\underline{X}) \in \Bbbk[X_1,\dots,X_n] \, , \; n -
\partial_{\bold{X}}(P_n) \geq \ell \;\, \forall\, n
\in \N \,\Big\}  $$
where  $ \, \check{x}_j := h^{-1} x_j \, $  and  $ \,
\partial_{\bold{X}}(f) \, $  denotes the degree of a
polynomial or a series  $ f $  in the  $ X_j $  ($ j =
1, \dots, n $),  \; whence we get
  $$  \displaylines{
   \quad  {F_h}^{\!\!\times} = \bigcup_{\ell \in \N} h^{-\ell}
I^\ell = \Big\{\, f = \sum\nolimits_{\bold{d} \in \N^{n+1}} \!
P_n(\check{x}_1, \dots, \check{x}_n) \cdot h^n \;\Big|   \hfill  \cr
   \hfill   \Big|\; P_n(X_1, \dots, X_n) \in \Bbbk[X_1,\dots,X_n] \, ,
\; \exists \, \ell_f \in \N : n - \partial_{\bold{X}}(P_n) \geq
\ell_f \; \forall\, n \in \N \,\Big\} \, .  \cr }  $$
Moreover, we easily see that  $ \, \bigcap_{\ell \in \N} h^\ell
{F_h}^{\!\times} = \{0\} \, $,  \, hence the natural completion
map  $ \, {F_h}^{\!\times} \longrightarrow {F_h}^{\!\vee} \, $
is an embedding.  Finally, when taking the  $ h $--adic
completion we get
  $$  {F_h}^{\!\!\vee} = \Big\{\, f = \sum\nolimits_{\bold{d} \in
\N^{n+1}} P_n(\check{x}_1,\dots,\check{x}_n) \cdot h^n \;\Big|\;
P_n(X_1,\dots,X_n) \in \Bbbk[X_1,\dots,X_n] \; \forall\, n \in \N
\,\Big\}  $$
that is  $ \, \displaystyle{ {F_h}^{\!\!\vee} = \Bbbk[\check{x}_1, \dots,
\check{x}_n][[h]] } \, $  as topological  $ \Bbbk[[h]] $--modules.

\vskip7pt

\proclaim{Proposition 3.2}  If  $ \, F_h[[\gerg]] \in \QFSHA \, $,
\, then  $ \, {F_h[[\gerg]]}^\vee \in \QUEA \, $.  Namely, we have
$ \; {F_h[[\gerg]]}^\vee = U_h(\gerg^*) \; $  (where  $ \gerg^* $
is the dual Lie bialgebra to  $ \gerg $),  that is
  $$  {F_h[[\gerg]]}^\vee \Big/ h \, {F_h[[\gerg]]}^\vee = U(\gerg^*)
\, .  $$
\endproclaim

\demo{Proof} \, Let  $ \, F_h[[\gerg]] \in \QFSHA \, $;  \, set for
simplicity  $ \, F_h := F_h[[\gerg]] \, $,  $ \, F_0 := F_h \big/ h \,
F_h = F[[\gerg]] \, $,  $ \, {F_h}^{\!\!\vee} := {F_h [[\gerg]]}^\vee $,
$ \, {F_0}^{\!\!\vee} := {F_h}^{\!\!\vee} \big/ h \, {F_h}^{\!\!\vee} $,
\, and let  $ \, \pi \colon \, F_h \relbar\joinrel\twoheadrightarrow
F_0 \, $  be the natural projection.
                                               \par
   From the discussion in \S 3.1, we recover the identification
$ \, F_h = \Bbbk[[x_1,\dots,x_n,h]] \, $  (for some  $ \, n \in \N \, $)
as topological  $ \Bbbk[[h]] $--modules,  where  $ \, x_j \in F_h \, $
for all  $ j $  and the  $ \, \bar{x}_j = \pi(x_j) \, $  gives
$ \, F_0 = \Bbbk[[\bar{x}_1,\dots,\bar{x}_n]] \, $  and generate  $ \,
\frak{m} := \text{\it Ker}\,(\epsilon_{F[[\gerg]]}) \, $.  Taking if
necessary  $ \, x_j - \epsilon(x_j) \, $  instead of  $ x_j $  (for
any  $ j \, $),  we can assume in addition that the  $ x_j $  belong
to  $ \, J := \text{\it Ker}\,(\epsilon_{\scriptscriptstyle F_h})
\, $,  \, so this kernel is the set of all formal series  $ f $
whose degree in the  $ x_j \, $,  \, call it  $ \partial_{\bold{x}}
(f) $,  is positive.  From \S 3.1 we also have  $ \, {F_h}^{\!\vee}
= \Bbbk[\check{x}_1,\dots,\check{x}_n][[h]] \, $  as topological
$ \Bbbk[[h]] $--modules.
                                               \par
   Since  $ F_0 $  is commutative, we have  $ \, x_i \, x_j - x_j \,
x_i = h \, \chi \, $  for some  $ \, \chi \in F_h \, $,  \, and in
addition we must have  $ \, \chi \in J \, $  too, thus  $ \, \chi =
\sum_{j=0}^n c_j(h) \cdot x_j + f(x_1,\dots,x_n,h) \, $  where  $ \,
c_j(h) \in \Bbbk[[h]] \, $  for all  $ j $  and  $ \, f(x_1,\dots,x_n,h)
\in \Bbbk[[x_1,\dots,x_n,h]] \, $  with  $ \, \partial_{\bold{x}}(f) > 1
\, $.  Then
  $$  \check{x}_i \, \check{x}_j - \check{x}_j \, \check{x}_i
= h^{-2} \cdot h \, \chi = \sum\nolimits_{j=0}^n c_j(h) \cdot
\check{x}_j + h^{-1} \check{f}(\check{x}_1,\dots,\check{x}_n,h)
\, ,  $$
where  $ \, \check{f}(\check{x}_1,\dots,\check{x}_n,h) \in
k[\check{x}_1,\dots,\check{x}_n][[h]] \, $  is formally obtained
from  $ f(x_1,\dots,x_n,h) $  simply by rewriting  $ \, x_j = h \,
\check{x}_j \, $  for all  $ j $.  Then since  $ \, \partial_{\bold{x}}
(f) > 1 \, $  we have  $ \, h^{-1} \check{f}(\check{x}_1, \dots,
\check{x}_n, h) \in h \, \Bbbk[\check{x}_1,\dots,\check{x}_n][[h]] \, $,
\; whence
  $$  \check{x}_i \, \check{x}_j - \check{x}_j \, \check{x}_i
\equiv \sum\nolimits_{j=0}^n c_j(h) \cdot \check{x}_j \, \mod\,
h \, {F_h}^{\!\!\vee} \, .  $$
This shows that  {\sl the  $ \Bbbk $--span  of the set of cosets  $ \,
{\big\{\, \check{x}_j \mod\, h \, {F_h}^{\!\!\vee} \,\big\}}_{j=1,
\dots, n} \,$  is a Lie algebra},  which we call  $ \gerh $.  Then
the identification  $ \, {F_h}^{\!\!\vee} = \Bbbk[\check{x}_1, \dots,
\check{x}_n][[h]] \, $  shows that  $ \, {F_0}^{\!\!\vee} = U(\gerh)
\, $,  \, so that  $ \, {F_h}^{\!\!\vee} \in \QUEA \, $,  \, q.e.d.
                                             \par
   Our purpose now is to prove that  $ \, \gerh \cong \gerg^* \, $
as Lie bialgebras.  For this we have to improve a bit the previous
analysis.  Recall that\footnote{ Hereafter, the product of ideals in
                     a topological algebra will be understood as the
                     {\sl closure\/}  of their algebraic product. }
       $ \, \gerg := {\big( \frak{m} \big/ \frak{m}^2 \big)}^* \, $,
\, that  $ \frak{m} $  (the unique maximal ideal of $ F[[\gerg]] \, $)
is closed under the Poisson bracket of  $ F[[\gerg]] $,  and that the
dual Lie bialgebra  $ \gerg^* $  can be realized as  $ \, \gerg^* =
\frak{m} \big/ \frak{m}^2 \, $,  \, its Lie bracket being induced
by the Poisson bracket.
                                             \par
   Consider  $ \, J^\vee := h^{-1} J \subset {F_h}^{\!\!\vee} \, $.
Multiplication by  $ h^{-1} $  yields a  $ \Bbbk[[h]] $--module
isomorphism  $ \, \mu \, \colon \, J \,{\buildrel \cong \over
{\lhook\joinrel\relbar\joinrel\twoheadrightarrow}}\, J^\vee $.
Furthermore, the specialisation map  $ \, \pi^\vee \colon \,
{F_h}^{\!\!\vee} \! \relbar\joinrel\twoheadrightarrow {F_0}^{\!\!\vee}
\! = \! U(\gerh) \, $  restricts to a similar map  $ \, \eta \, \colon
J^\vee \! \relbar\joinrel\twoheadrightarrow \Cal{J}_{\,0}^\vee \! :=
\! J^\vee \Big/ J^\vee \cap {\big( h \, {F_h}^{\!\!\vee} \big)} \, $.
The latter has kernel  $ \, J^\vee \cap {\big( h \, {F_h}^{\!\!\vee}
\big)} \, $:  \, we contend that this is equal to  $ \, \big( J +
J^\vee J \, \big) \, $.  In fact, let $ \, y \in J^\vee \cap h \,
{F_h}^{\!\!\vee} \, $:  then the series  $ \, \gamma := h \, y \in
J \, $  has  $ \, \partial_{\bold{x}}(\gamma) > 0 \, $.  As above we
write  $ \, y = h^{-1} \gamma \, $  as  $ \, y = h^{-1} \gamma \in
{F_h}^{\!\!\vee} = \Bbbk[\check{x}_1,\dots,\check{x}_n][[h]] \, $:  \,
then  $ \, y = h^{-1} \gamma \in h \, \Bbbk[\check{x}_1,\dots,\check{x}_n]
[[h]] \, $  means  $ \, \partial_{\bold{x}}(\gamma) > 1 \, $,  \,
or  $ \, \partial_{\bold{x}}(\gamma) = 1 \, $  {\sl and\/}  $ \,
\partial_h(\gamma) > 0 \, $  (i.e.~$ \, \gamma \in h \, \Bbbk[[x_1,
\dots,x_n,h]] \, $),  i.e.~exactly  $ \, \gamma \in h \, J + J^2
\, $,  \, so  $ \, y \in J + J^\vee J $,  \, which proves our claim
true.  Note also that  $ \, \eta \big( J^\vee \big) = \eta \big(
\oplus_{j=1}^n \! \Bbbk[[h]] \cdot x_j \big) = \oplus_{j=1}^n \Bbbk \cdot
\check{x}_j = \gerh \, $.
                                             \par
   Now, recall that  $ \, \gerg^* = \frak{m} \big/ \frak{m}^2 \, $:
\, we fix a  $ \Bbbk $--linear  section  $ \; \nu \, \colon \, \gerg^*
\lhook\joinrel\longrightarrow \frak{m} \, $  of the projection
$ \; \rho \, \colon \, \frak{m} \relbar\joinrel\twoheadrightarrow
\frak{m} \big/ \frak{m}^2 = \gerg^* \, $  such that  $ \, \gamma
(\frak{m}^2) \subseteq h \, J + J^2 \, $.  Moreover, the specialisation
map  $ \; \pi \, \colon \, F_h \relbar\joinrel\twoheadrightarrow
F_0 \; $  restricts to  $ \; \pi' \colon \, J
\relbar\joinrel\twoheadrightarrow J \big/ (J \cap
h \, F_h) = J \big/ h \, J = \frak{m} \, $;  \, we fix a
$ \Bbbk $--linear section  $ \, \gamma \, \colon \, \frak{m}
\lhook\joinrel\relbar\joinrel\rightarrow J \, $  of  $ \pi' $.
Now consider the composition map  $ \, \sigma := \eta \circ \mu \circ
\gamma \circ \nu \, \colon \, \gerg^* \longrightarrow \gerh \, $.
This is well-defined, i.e.~it is independent of the choice of  $ \nu $
and  $ \gamma \, $.  Indeed, if  $ \, \nu, \nu' \colon \, \gerg^*
\lhook\joinrel\relbar\joinrel\rightarrow J_{\scriptscriptstyle F} \, $
are two sections of  $ \rho $,  and  $ \sigma $,  $ \sigma' $  are
defined correspondingly (with the same fixed  $ \gamma $  for both),
then  $ \, \hbox{\it Im}\,(\nu-\nu') \subseteq \hbox{\it Ker}\,(\rho)
= {J_{\scriptscriptstyle F}}^{\!2} \subseteq \hbox{\it Ker}\, (\eta
\circ \mu \circ \gamma) \, $,  \, so that  $ \, \sigma = \eta \circ
\mu \circ \gamma \circ \nu = \eta \circ \mu \circ \gamma \circ \nu'
= \sigma' \, $.  Similarly, if  $ \, \gamma, \gamma' \colon \,
J_{\scriptscriptstyle F} \lhook\joinrel\relbar\joinrel\rightarrow
J \, $  are two sections of  $ \pi' $,  and  $ \sigma $,  $ \sigma' $
are defined correspondingly (with the same  $ \nu $  for both), we
have  $ \, \hbox{\it Im}\,(\gamma-\gamma') \subseteq \hbox{\it Ker}\,
(\pi') = h \, J \subseteq h \, J + J^2 = \hbox{\it Ker}\,(\eta \circ
\mu) \, $,  thus  $ \, \sigma = \eta \circ \mu \circ \gamma \circ \nu
= \eta \circ \mu \circ \gamma' \circ \nu = \sigma' \, $,  \, q.e.d.
In a nutshell,  $ \sigma $  is the composition map
  $$  \gerg^*  \hskip7pt  {\buildrel \bar{\nu} \over
{\lhook\joinrel\relbar\joinrel\relbar\joinrel\relbar%
\joinrel\twoheadrightarrow}}
\hskip7pt  J_{\scriptscriptstyle F} \Big/
{J_{\scriptscriptstyle F}}^{\! 2}  \hskip7pt
{\buildrel \bar{\gamma} \over
{\lhook\joinrel\relbar\joinrel\relbar\joinrel\relbar%
\joinrel\twoheadrightarrow}}
\hskip7pt  J \Big/ \! \big( J^2 + h \hskip1pt J \, \big)
\hskip7pt {\buildrel \bar{\mu} \over
{\lhook\joinrel\relbar\joinrel\relbar\joinrel\relbar%
\joinrel\twoheadrightarrow}}
\hskip7pt  J^\vee \Big/ \big( J + J^\vee J \big)
\hskip7pt  {\buildrel \bar{\eta} \over
{\lhook\joinrel\relbar\joinrel\relbar\joinrel\relbar%
\joinrel\twoheadrightarrow}}
\hskip7pt  \gerh  $$
where the maps  $ \bar{\nu} $,  $ \bar{\gamma} $,  $ \bar{\mu} $,
$ \bar{\eta} $,  are the ones canonically induced by  $ \nu $,
$ \gamma $,  $ \mu $,  $ \eta $,  and  $ \bar{\nu} $,  {\sl
resp.~$ \bar{\gamma} $,  does not depend on the choice of  $ \nu $,
resp.~$ \gamma $},  as it is the inverse of the isomorphism
$ \, \bar{\rho} \, \colon \, J_{\scriptscriptstyle F} \big/
{J_{\scriptscriptstyle F}}^{\! 2} \,{\buildrel \cong \over
{\lhook\joinrel\relbar\joinrel\twoheadrightarrow}}\, \gerg^* \, $,
resp.~$ \, \overline{\pi'} \, \colon \, J \Big/ \big( J^2 + h \, J \,
\big) \,{\buildrel \cong \over
{\lhook\joinrel\relbar\joinrel\twoheadrightarrow}}\,
J_{\scriptscriptstyle F} \big/ {J_{\scriptscriptstyle F}}^{\!2}
\, $,  induced by  $ \, \rho \, $,  resp.~by  $ \pi' $.  We use
this remark to show that  $ \sigma $  is also an isomorphism of
the Lie bialgebra structure.
                                             \par

   Using the vector space isomorphism  $ \, \sigma \, \colon \,
\gerg^* \,{\buildrel \cong \over \longrightarrow}\, \gerh \, $
we pull-back the Lie bialgebra structure of  $ \gerh $  onto
$ \gerg^* $,  and denote it by  $ \big( \gerg^*, {[\,\ ,\ ]}_\bullet,
\delta_\bullet \big) $;  \, on the other hand,  $ \gerg^* $  also
carries its natural structure of Lie bialgebra, dual to that of
$ \gerg $  (e.g., the Lie bracket is induced by restriction of
$ \{ \hskip4pt , \hskip3pt \} \, $),  denoted by  $ \big( \gerg^*,
{[\,\ ,\ ]}_\times, \delta_\times \big) $:  we must prove that these
two structures coincide.
                                             \par
   First,  {\it for all  $ \, x_1 $,  $ x_2 \in \gerg^* \, $  we have
$ \; {\big[x_1,x_2\big]}_\bullet = {\big[x_1,x_2\big]}_\times \, $}.
                                             \par
   Indeed, let  $ \, f_i := \nu(x_i) \, $,  $ \, \varphi_i := \gamma(f_i)
\, $,  $ \, \varphi^\vee_i := \mu(\varphi_i) \, $,  $ \, y_i := \eta
\big( \varphi^\vee_i \big) \, $  \, ($ \, i = 1, 2 $).  Then
  $$  \displaylines{
   {\big[ x_1, x_2 \big]}_\bullet := \sigma^{-1} \Big( {\big[
\sigma(x_1), \sigma(x_2) \big]}_\gerh \Big) = \sigma^{-1} \big(
[y_1,y_2] \big)  = \big( \rho \circ \pi' \circ \mu^{-1} \big) \Big(
\big[ \varphi^\vee_1, \varphi^\vee_2 \big] \Big) =   \hfill   \cr
    \hfill   = \big( \rho \circ \pi' \big) \big( h^{-1}
[\varphi_1, \varphi_2] \big) = \rho \big( \{f_1,f_2\} \big)
=: {\big[ x_1, x_2 \big]}_\times \; ,  \qquad  \text{q.e.d.}
\cr }  $$
   \indent   The case of cobrackets can be treated similarly; but
since they take values in tensor squares, we make use of suitable
maps  $ \, \nu_\otimes := \nu^{\otimes 2} $,  $ \, \gamma_\otimes
:= \gamma^{\otimes 2} $,  etc.; we set also  $ \, \chi_\otimes \!
:= \eta_\otimes \! \circ \mu_\otimes = {(\eta \circ \mu)}^{\otimes 2}
\, $  and  $ \, \nabla \! := \Delta - \Delta^{\text{op}} $.  Then
{\it for all  $ \, x \in \gerg^* \, $  we have  $ \; \delta_\bullet(x)
= \delta_\times(x) \, $}.
                                             \par
   Indeed, let  $ \, f := \nu(x) \, $,  $ \, \varphi := \gamma(f) \, $,
$ \, \varphi^\vee := \mu(\varphi) \, $,  $ \, y := \eta \big(
\varphi^\vee \big) \, $.  Then we have
  $$  \displaylines{
   \delta_\bullet(x) := {\sigma_\otimes}^{\!-1} \big( \delta_\gerh
(\sigma(x)) \big) = {\sigma_\otimes}^{\!\!-1} \big( \delta_\gerh
\big( \eta \big( \varphi^\vee \big) \big) \big) =
{\sigma_\otimes}^{\!\!-1} \Big( \eta_\otimes \big(
h^{-1} \nabla \big( \varphi^\vee \big) \big) \Big) =   \hfill  \cr
   = {\sigma_\otimes}^{\!\!-1} \Big( {\big( \eta \circ \mu
\big)}_\otimes \big( \nabla(\varphi) \big) \Big) =
{\sigma_\otimes}^{\!\!-1} \Big( {\big( \eta \circ \mu \circ \gamma
\big)}_\otimes \big( \nabla(f) \big) \Big) = \rho \big( \nabla(f)
\big) = \rho \big( \nabla(\nu(x)) \big) = \delta_\times(x)
\cr }  $$
where the last equality holds because  $ \delta_\times(x) $  is
uniquely defined as the unique element in  $ \, \gerg^* \otimes
\gerg^* \, $  such that  $ \; \big\langle u_1 \otimes u_2 \, , \,
\delta_\times(x) \big\rangle = \big\langle [u_1,u_2] \, , \, x
\big\rangle \; $  for all  $ \, u_1, u_2 \! \in \! \gerg \, $,
\, and we have
  $$  \big\langle [u_1,u_2] \, , \, x \big\rangle = \big\langle
[u_1,u_2] \, , \, \rho(f) \big\rangle = \big\langle u_1 \otimes u_2
\, , \nabla(f) \big\rangle = \big\langle u_1 \otimes u_2 \, ,
\rho\big(\nabla(\nu(x))\big) \big\rangle \; .   \quad  \square  $$
\enddemo

   Now we need one more technical lemma.  From now on, if  $ \gerg $
is any Lie algebra and  $ \, \bar{x} \in U(\gerg) \, $,  we denote by
$ \, \partial \big( \bar{x} \big) \, $  the degree of  $ \bar{x} $
w.r.t.~the standard filtration of  $ U(\gerg) $.

\vskip7pt

\proclaim{Lemma 3.3} Let  $ U_h $  be a QUEA, let  $ \, x' \in
{U_h}^{\!\prime} \, $,  and let  $ \, x \in U_h \setminus h \, U_h
\, $,  $ \, n \in \N \, $,  be such that  $ \, x' = h^n x \, $.  Set
$ \, \bar{x} := x \,\mod h \, U_h \in {\big(U_h\big)}_0 \, $.
Then  $ \, \partial(\bar{x}) \leq n \, $.
\endproclaim

\demo{Proof}  (cf.~[EK], Lemma 4.12)  By hypothesis  $ \, \delta_{n+1}
(x') \in h^{n+1} {U_h}^{\! \otimes (n+1)} $,  hence  $ \, \delta_{n+1}
(x) \in h \, {U_h}^{\! \otimes (n+1)} \, $,  therefore  $ \,
\delta_{n+1}(\bar{x}) = 0 \, $,  i.e.~$ \, \bar{x} \in Ker \big(
\delta_{n+1} \colon \, U(\gerg) \longrightarrow {U(\gerg)}^{\otimes
(n+1)} \big) \, $,  where  $ \gerg $  is the Lie bialgebra such that
$ \, {\big( U_h \big)}_0 := U_h \big/ h \, U_h = U(\gerg) \, $.  But
the latter kernel equals the subspace  $ \, {U(\gerg)}_n := \big\{\,
\bar{y} \in U(\gerg) \,\big\vert\, \partial(\bar{y}) \leq n \,\big\}
\, $  (cf.~[KT], \S 3.8), whence the claim follows.   $ \square $
\enddemo

\vskip7pt

\proclaim{Proposition 3.4} Let  $ F_h $  be a QFSHA.
Then  $ \, {\big( {F_h}^{\!\vee} \big)}' = F_h \, $.
\endproclaim

\demo{Proof}  As a matter of notation, we set  $ \,
{J_{\scriptscriptstyle F_h}}^{\!\!\vee} := h^{-1}
J_{\scriptscriptstyle F_h} \, $,  \, and we denote by
$ \, \bar{x} \in {F_0}^{\!\vee} \, $  the image of any
$ \, x \in {F_h}^{\!\vee} \, $  inside  $ \, {F_h}^{\!\vee}
\big/ h \, {F_h}^{\!\vee} = {F_0}^{\!\vee} $.
                                                 \par
   Now, for any  $ \, n \in \N $,  \, we have  $ \, \delta_n(F_h)
\subseteq {J_{\scriptscriptstyle F_h}}^{\hskip-3pt \otimes n} \, $
(see \S 1.5); this can be read as  $ \, \delta_n(F_h) \subseteq
{J_{\scriptscriptstyle F_h}}^{\hskip-3pt \otimes n} = h^n {\big(
h^{-1} J_{\scriptscriptstyle F_h} \big)}^{\! \otimes n} \subseteq
h^n {\big( {F_h}^{\!\times} \big)}^{\! \otimes n} \subseteq h^n
{\big( {F_h}^{\!\vee} \big)}^{\! \otimes n} \, $,  which gives
$ \, F_h \subseteq {\big( {F_h}^{\!\vee} \big)}' \, $.
                                                 \par
   Conversely, let  $ \, x' \in {\big( {F_h}^{\!\vee} \big)}' \setminus
\{0\} \, $  be given; as  $ \, {F_h}^{\!\vee} \in {\Cal T}_{\otimeshat}
\, $,  \, there are (unique)  $ \, n \in \N $,  $ \, x \in {F_h}^{\!
\vee} \setminus h \, {F_h}^{\!\vee} $,  \, such that  $ \, x' = h^n x
\, $.  By Proposition 3.2,  $ {F_h}^{\!\vee} $  is a QUEA, with
semiclassical limit  $ U(\gerh) $  where  $ \, \gerh = \gerg^*
\, $ if  $ \, F_0 = F[[\gerg]] \, $.  Fix an ordered basis  $ \,
{\{b_\lambda\}}_{\lambda \in \Lambda} \, $  of  $ \gerh $  and
a subset  $ \, {\{x_\lambda\}}_{\lambda \in \Lambda} \, $  of
$ {F_h}^{\!\vee} $  such that  $ \, \overline{x}_\lambda = b_\lambda
\, $  for all  $ \, \lambda \, $;  \, in particular, since  $ \, \gerh
\subset \hbox{\it Ker}\,(\epsilon_{\scriptscriptstyle U(\gerg)}) \, $
we can choose the  $ x_\lambda $  inside  $ \, {J_{\scriptscriptstyle
F_h}}^{\!\!\vee} := h^{-1} J_{\scriptscriptstyle F_h} \, $:  \, so
$ \, x_\lambda = h^{-1} x'_\lambda \, $  for some  $ \, x'_\lambda
\in J_{\scriptscriptstyle F_h} \, $,  \, for all  $ \, \lambda
\in \Lambda \, $.
                                             \par
   Since Lemma 3.3 gives  $ \, \partial(\bar{x}) \leq n \, $,  \, that
is  $ \, \bar{x} \in {U(\gerg)}_n := \big\{\, \bar{y} \in U(\gerg)
\,\big\vert\, \partial(\bar{y}) \leq n \,\big\} \, $,  \, by the
PBW theorem we can write  $ \bar{x} $  as a polynomial  $ \, P
\big( {\{b_\lambda\}}_{\lambda \in \Lambda} \big) \, $  in the
$ b_\lambda $  of degree  $ \, d \leq n \, $  (with coefficients
in  $ \Bbbk \, $);  then  $ \, x_0 := P \big( {\{x_\lambda\}}_{\lambda
\in \Lambda} \big) \equiv x \mod\, h \, {F_h}^{\!\vee} \, $,  \, that
is  $ \, x = P \big( {\{x_\lambda\}}_{\lambda \in \Lambda} \big)
+ h \, x_{\langle 1  \rangle} \, $  for some  $ \, x_{\langle 1
\rangle} \in H^\vee \, $.  Now we can write  $ \, x_0 := P \big(
{\{x_\lambda\}}_{\lambda \in \Lambda} \big) = \sum_{s=0}^d h^{-s}
j_s \; ( \, \in \! {F_h}^{\!\times} ) $,  \, where every  $ \, j_s
\in {J_{\scriptscriptstyle F_h}}^{\!s} \, $  is a homogeneous
polynomial in the  $ x'_\lambda $  of degree  $ s $,  and  $ \,
j_d \not= 0 \, $;  but then  $ \, h^n x_0 = \sum_{s=0}^d h^{n-s}
j_s \in F_h \, $  because  $ \, d \leq n \, $.  Since  $ \, F_h
\subseteq {\big( {F_h}^{\!\vee} \big)}' \, $   --- thanks to the
first part of the proof ---   we get also  $ \, h^{n+1} x_{\langle
1 \rangle} = h^n (x - x_0) = x' - h^n x_0 \in {\big( {F_h}^{\!\vee}
\big)}' \, $:  thus
  $$  x' = h^n x_0 + h^{n+1} x_{\langle 1 \rangle} \, ,  \qquad
\hbox{with}  \quad  h^n x_0 \in F_h  \quad  \hbox{and}  \quad
h^{n+1} x_{\langle 1 \rangle} \in {\big( {F_h}^{\!\vee} \big)}'
\, .  $$
   \indent   If  $ \, x_{(1)} := h^{n+1} x_{\langle 1 \rangle}
\, $  is zero we are done; if not, we can repeat the argument for
$ x_{(1)} $  in the role of  $ \, x_{(0)} := x' \, $:  this will
provide us with an  $ \, x_1 \in {F_h}^{\!\times} \, $  and an
$ \, x_{\langle 2 \rangle} \in {F_h}^{\!\vee} \, $  such that
$ \, x_{(1)} = h^{n+1} x_1 + h^{n+2} x_{\langle 2 \rangle} \, $,
\, with  $ \, h^{n+1} x_1 \in F_h \, $  and  $ \, h^{n+2} x_{\langle
2 \rangle} \in {\big( {F_h}^{\!\vee} \big)}' \, $.  Iterating, we
eventually find a sequence  $ \, {\{x_\ell\}}_{\ell \in \N} \subset
{F_h}^{\!\times} \, $  such that  $ \, h^{n+\ell} x_\ell \in
F_h \, $  for all  $ \, \ell \in \N \, $,  and  $ \, x' =
\sum\limits_{\ell=0}^{+\infty} h^{n+\ell} x_\ell \, $,  \,
in the sense that the right-hand-side series does converge
to  $ x' $  inside  $ {F_h}^{\!\vee} $.  Furthermore, this
convergence takes place inside  $ F_h $  as well: indeed, the very
construction gives  $ \, h^{n+\ell} x_\ell = h^{n+\ell} P_{d_\ell}
\big( {\{x_\lambda\}}_{\lambda \in \Lambda} \big) = h^{n+\ell}
P_{d_\ell} \big( {\{h^{-1} x'_\lambda\}}_{\lambda \in \Lambda}
\big) \, $  (where  $ P_{d_\ell} $  is a suitable polynomial of
degree  $ \, d_\ell \leq n + \ell \, $)  and this last element
belongs to  $ \, {I_{\scriptscriptstyle F_h}}^{\!\! n + \ell} \, $:
\, but  $ F_h $  is a QFSHA, hence it is complete w.r.t.~the
$ I_{\scriptscriptstyle F_h} \! $--adic  topology, so the series
$ \, x' = \sum\limits_{\ell=0}^{+\infty} h^{n+\ell} x_\ell \, $
does converge (to  $ x' $)  inside  $ F_h $,  q.e.d.   $ \square $
\enddemo

\vskip7pt

   {\bf 3.5 An explicit description of  $ \, {U_h(\gerg)}' \, $  (for
\, \hbox{\bfit char}($\Bbbk$) = 0 ).}  \, When  $ \, \text{\it char}\,
(\Bbbk) = 0 \, $,  \, for any  $ \, U_h(\gerg) \in \QUEA \, $  we can give
an explicit description of  $ \, {U_h(\gerg)}' \, $,  \, as follows.
                                                     \par
   Like in the proof of Proposition 2.6 consider  $ \, F_h[[\gerg]]
:= {U_h(\gerg)}^* \in \HA_{\otimestilde} \, $  and its natural Hopf
pairing with  $ U_h(\gerg) $: then we showed that  $ \, {U_h(\gerg)}'
= \Big( \! {\big( F_h[[\gerg]] \big)}^{\!\vee} \Big)^{\!*} \, $.  Note
that this time we have in addition  $ \, F_h[[\gerg]] \in \QFSHA \, $,
\, with  $ \, F_0[[\gerg]] := F_h[[\gerg]] \big/ h \, F_h[[\gerg]] =
{U(\gerg)}^* = F[[\gerg]] \, $.
                                                     \par
   Pick any basis  $ \, {\{ \overline{x}_i \}}_{i \in \Cal{I}} \, $
of  $ \gerg $,  endowed with some total order; then (PBW theorem) the
set of ordered monomials  $ \, {\big\{ \overline{x}^{\,\underline{e}}
\big\}}_{\underline{e} \in {( \N^\Cal{I} )}_0} \, $  is a basis of
$ U(\gerg) $:  hereafter,  $ {\big( \N^\Cal{I} \,\big)}_0 $  denotes
the set of functions from  $ \Cal{I} $  to  $ \N $  with finite
support, and  $ \, \overline{x}^{\,\underline{e}} := \prod_{i \in
\Cal{I}} \overline{x}_i^{\,\underline{e}(i)} \, $  for all
$ \, \underline{e} \in {\big( \N^\Cal{I} \,\big)}_0 \, $
and all indeterminates  $ \, x_1, \dots, x_n \, $.
Let  $ \, {\{ \overline{y}_i\}}_{i \in \Cal{I}} \, $  be the
          pseudobasis\footnote{ From now on, this means that each
element of  $ \gerg^* $  can be written uniquely as a (possibly
infinite) linear combination of elements of the pseudobasis: such
a (possibly infinite) sum will be convergent in the weak topology
of  $ \gerg^* $,  so a pseudobasis is a  {\sl topological basis\/}
w.r.t.~the weak topology. }
of  $ \gerg^* $  dual to  $ \, {\{ \overline{x}_i \}}_{i \in \Cal{I}}
\, $,  endowed with the same total order; then the set of "rescaled"
ordered monomials  $ \, {\big\{ \overline{y}^{\,\underline{e}} \big/
\underline{e}\,! \big\}}_{\underline{e} \in {( \N^\Cal{I} )}_0} \, $
(with  $ \, \underline{e}\,! := \prod_{i \in \Cal{I}} \underline{e}(i)!
\, $;  \,  that's where we need  $ \text{\it char}\,(\Bbbk) = 0 \, $)  \,
is the pseudobasis of  $ \, {U(\gerg)}^* = F_h[[\gerg]] \, $  dual
to the PBW basis  $ \, {\big\{ \overline{x}^{\,\underline{e}}
\big\}}_{\underline{e} \in {( \N^\Cal{I} )}_0} \, $  of  $ \, U(\gerg)
\, $,  \, namely  $ \, \big\langle \overline{x}^{\,\underline{e}},
\overline{y}^{\,\underline{e}'} \big/ \underline{e}\,! \big\rangle =
\delta_{\underline{e}, \underline{e}'} \, $  for all  $ \, \underline{e}
\, $,  $ \underline{e}' \in {\big( \N^\Cal{I} \big)}_0 \, $.
                                                     \par
   Lift  $ \, {\{\overline{x}_i\}}_{i \in \Cal{I}} \, $  to a subset
$ \, {\{x_i\}}_{i \in \Cal{I}} \subseteq U_h(\gerg) \, $  such that
$ \, \overline{x}_i = x_i \mod h \, U_h(\gerg) \, $,  \, and
$ \, {\{\overline{y}_i\}}_{i \in \Cal{I}} \, $  to a subset  $ \,
{\{y_i\}}_{i \in \Cal{I}} \subseteq \! F_h[[\gerg]] \, $  such
that  $ \, \overline{y}_i = y_i \mod h \, F_h[[\gerg]] \, $:  \,
then  $ \, {\big\{ x^{\,\underline{e}} \big\}}_{\underline{e} \in
{(\N^\Cal{I})}_0} \, $  is a topological basis of  $ U_h(\gerg) $
(as a topological  $ \Bbbk[[h]] $--module)  and similarly  $ \, {\big\{
y^{\,\underline{e}} \big/ \underline{e}\,! \big\}}_{\underline{e} \in
{(\N^\Cal{I})}_0} $  is a topological pseudobasis of  $ F_h[[\gerg]] $
(as a topological  $ \Bbbk[[h]] $--module),  and they are dual to
each other modulo  $ \, h \, $,  \, i.e.~$ \, \big\langle
x^{\,\underline{e}}, y^{\,\underline{e}'} \big/ \underline{e}\,!
\big\rangle \in \delta_{\underline{e},\underline{e}'} + h \, \Bbbk[[h]]
\, $  for all  $ \, \underline{e} \, $,  $ \underline{e}' \in {\big(
\N^\Cal{I} \big)}_0 \, $.  In addition,  $ \, y^{\,\underline{e}'} \big/
\underline{e}\,! \in {\big( I_{F_h[[\gerg]]} \big)}^{|\underline{e}'|}
\, $  for all  $ \, \underline{e}' \in {\big(\N^\Cal{I}\big)}_0 \, $,
\, where  $ \, \vert \underline{e}' \vert := \sum_{i \in \Cal{I}}
\underline{e}'(i) \, $.  Now,  $ U_h(\gerg) $  also contains a
topological basis dual to  $ \, {\big\{ y^{\,\underline{e}} \big/
\underline{e}\,! \big\}}_{\underline{e} \in {(\N^\Cal{I})}_0} \, $,
\, call it  $ \, {\big\{\eta_{\,\underline{e}} \,\big\}}_{\underline{e}
\in {( \N^\Cal{I} )}_0} \, $:  \, indeed, from the previous analysis
we see   --- by the "duality mod  $ h $"  mentioned above ---
that such a basis is given by  $ \; \eta_{\underline{e}} \, =
\, x^{\,\underline{e}} \, + \, \sum_{n=1}^{+\infty} h^n \! \cdot
\sum_{\underline{e}' \in \, {( \N^\Cal{I} )}_0} c^{(n)}_{\underline{e},
\underline{e}'}  \hskip2pt  x^{\,\underline{e}'} \, $  for some
$ \, c^{(n)}_{\underline{e},\underline{e}'} \in \Bbbk \; $  ($ \,
\underline{e} $,  $ \underline{e}' \in {\big( \N^\Cal{I}
\big)}_0 \, $),  \, so  $ \, {\big\{ \eta_{\,\underline{e}}
\,\big\}}_{\underline{e} \in {(\N^\Cal{I})}_0} \, $  is a lift
of the PBW basis  $ \, {\big\{ \overline{x}^{\, \underline{e}}
\big\}}_{\underline{e} \in {( \N^\Cal{I} )}_0} \, $  of  $ \,
U(\gerg) \, $.  Since  $ \, {\big\{ y^{\,\underline{e}} \big/
\underline{e}\,! \big\}}_{\underline{e} \in {( \N^\Cal{I} )}_0} \, $
is a topological pseudobasis of  $ \, F_h[[\gerg]] \, $  and  $ \,
y^{\, \underline{e}'} \big/ \underline{e}\,! \in {\big( J_{F_h[[\gerg]]}
\big)}^{\vert \underline{e}' \vert} \, $  for all  $ \underline{e}' $,
the set  $ \, {\big\{ h^{-|\underline{e}|} \, y^{\,\underline{e}} \big/
\underline{e}\,! \big\}}_{\underline{e} \in {(\N^\Cal{I})}_0} \, $  is
a topological basis of the topologically free  $ \Bbbk[[h]] $--module
$ {\big(F_h[[\gerg]]\big)}^\vee \, $:  \, then the dual pseudobasis of
$ \, \Big( \! {\big( F_h[[\gerg]] \big)}^\vee \Big)^* = {U_h(\gerg)}'
\, $  to this basis is  $ \, {\big\{ h^{\vert \underline{e} \vert} \,
\eta_{\,\underline{e}} \,\big\}}_{\underline{e} \in {(\N^\Cal{I})}_0}
\, $,  \, so  $ \, {U_h(\gerg)}' \, $  is the set  $ \; \Big\{\,
\sum_{\underline{e} \in {( \N^\Cal{I} )}_0}  \hskip-1pt
a_{\underline{e}} \, h^{|\underline{e}|} \, \eta_{\underline{e}}
\hskip7pt \Big| \hskip5pt a_{\underline{e}} \in \Bbbk[[h]] \;\; \forall
\;\; \underline{e} \in {( \N^\Cal{I} )}_0 \,\Big\} \, $.
                                                     \par
   Now observe that  $ \; {(h \, x)}^{\underline{e}} =
h^{|\underline{e}|} \, x^{\underline{e}} \, \equiv
h^{|\underline{e}|} \, \eta_{\underline{e}} \mod h
\, {U_h(\gerg)}' \, $  by construction; therefore
  $$  \Big\{ \sum\nolimits_{\underline{e} \in {(\N^\Cal{I})}_0}
\hskip-3pt  a_{\underline{e}} \, {(h \, x)}^{\underline{e}}
\hskip3,8pt  \Big|  \hskip3pt  a_{\underline{e}} \in \Bbbk[[h]] \;\,
\forall \; \underline{e} \,\Big\} =  \hskip-1,5pt  \Big\{
\sum\nolimits_{\underline{e} \in {(\N^\Cal{I})}_0}  \hskip-3pt
a_{\underline{e}} \, h^{|\underline{e}|} \, \eta_{\underline{e}}
\hskip3,8pt  \Big|  \hskip3pt  a_{\underline{e}} \in \Bbbk[[h]] \;\,
\forall \; \underline{e} \,\Big\}  \hskip-1pt  =
\hskip-1pt  {U_h(\gerg)}' \, .  $$
Finally, up to taking  $ \, x_i - \epsilon(x_i) \, $,  \, one can
also choose the  $ \, x_i \, $  so that  $ \, \epsilon(x_i) = 0 \, $.
                                                     \par
   To summarize, the outcome is the following:

\vskip3pt

   {\it  Given any basis of  $ \gerg $,  there exists a lift  $ \,
{\{x_i\}}_{i \in \Cal{I}} \, $  of it in  $ U_h(\gerg) $  such that
$ \, \epsilon(x_i) = 0 \, $  and  $ {U_h(\gerg)}' $  is nothing but
the topological  $ \Bbbk[[h]] $--algebra in  $ \Cal{P}_{\otimestilde} $
generated (in topological sense) by  $ \, {\{h \, x_i\big\}}_{i
\in \Cal{I}} \, $,  \, thus  $ \, {U_h(\gerg)}' = \Big\{
\sum_{\underline{e} \in {(\N^\Cal{I})}_0} \hskip-1pt
a_{\underline{e}} \, h^{|\underline{e}|} \, x^{\underline{e}}
\hskip5pt \Big| \hskip3pt a_{\underline{e}} \in \Bbbk[[h]]
\;\, \forall \; \underline{e} \,\Big\} \, $  as a
subset of  $ \, U_h(\gerg) $.}

\vskip3pt

   $ \underline{\hbox{\sl Remark:}} $ \, this description of
$ {U_h(\gerg)}' $  implies that  {\sl the weak topology on}
$ {U_h(\gerg)}' $,  which coincides with its  $ I_{\scriptscriptstyle
{U_h(\gerg)}'} $--adic  topology,  {\sl does coincide with the
induced topology\/}  (of  $ {U_h(\gerg)}' $  as a subspace of
$ U_h(\gerg) $,  the latter being endowed with the  $ h $--adic
topology).  This defines the topology on  $ {U_h(\gerg)}'$  in an
intrinsic way, i.e.~without referring to any identification of
$ {U_h(\gerg)}' $  with the dual space to some  $ \, X \! \in
\Cal{T}_{\otimeshat} $  (as we did instead to prove Proposition 2.6).

\vskip7pt

\proclaim{Proposition 3.6}  Assume  $ \, \text{\it char}\,(\Bbbk)
= 0 \, $.  If  $ \, U_h(\gerg) \in \QUEA \, $,  \, then
$ \, {U_h(\gerg)}' \in \QFSHA \, $.  Namely, we have
$ \; {U_h(\gerg)}' = F_h[[\gerg^*]] \; $  (where  $ \gerg^* $
is the dual Lie bialgebra to  $ \gerg $), that is
  $$  {U_h(\gerg)}' \Big/ h \, {U_h(\gerg)}' = F[[\gerg^*]] \, .  $$
\endproclaim

\demo{Proof}  Consider  $ \, F_h[[\gerg]] := {U_h(\gerg)}^*
\! \in \QFSHA \, $  (cf.~Lemma 2.1); then  $ \, U_h(\gerg)
\hskip5pt {\buildrel \, {h \rightarrow 0} \, \over
{\relbar\joinrel\llongrightarrow}} \hskip4pt  U(\gerg) \, $
implies  $ \, F_h[[\gerg]] := {U_h(\gerg)}^* \hskip1pt {\buildrel
\, {h \rightarrow 0} \, \over {\relbar\joinrel\llongrightarrow}}
\hskip4pt {U(\gerg)}^* = F[[\gerg]] \, $.  By Proposition 3.2,
$ {F_h[[\gerg]]}^\vee $  is a QUEA, with semiclassical limit
$ U(\gerg^*) \, $;  \, by Proposition 2.2 we have  $ \, {U_h(\gerg)}'
= {\big( {F_h[[\gerg]]}^\vee \big)}^* $,  \; thus  $ {U_h(\gerg)}' $
is a QFSHA with semiclassical limit  $ F[[\gerg^*]] $:  \, indeed,
$ \, {F_h[[\gerg]]}^\vee \hskip1pt {\buildrel \, {h \rightarrow 0}
\, \over {\relbar\joinrel\llongrightarrow}} \hskip4pt  U(\gerg^*) \, $
implies  $ \, {U_h(\gerg)}' \! = {\big( {F_h[[\gerg]]}^\vee \big)}^*
\hskip0pt {\buildrel \, {h \rightarrow 0} \, \over
{\relbar\joinrel\loongrightarrow}} \hskip3pt {U(\gerg^*)}^*
\! = F[[\gerg^*]] \, $,  \, i.e.~$ \, {U_h(\gerg)}' \Big/ h \,
{U_h(\gerg)}' = F[[\gerg^*]] \, $  as  \hbox{claimed.   $ \square $}
\enddemo

\vskip7pt

\proclaim{Proposition 3.7}  Assume  $ \, \text{\it char}\,(\Bbbk)
= 0 \, $.  Let  $ U_h $  be a QUEA.  Then  $ \, {\big( {U_h}'
\,\big)}^{\!\vee} = U_h \, $.
\endproclaim

\demo{Proof}  Consider  $ \, F_h := {U_h\phantom{|}}^{\!\!\!*}
\in \QFSHA \, $  (cf.~Lemma 2.1): then Proposition 3.4 yields
$ \, {\big( {F_h}^{\!\!\vee} \big)}' = F_h \, $;  \, furthermore,
$ \, U_h = {\big( {U_h\phantom{|}}^{\!\!\!*} \,\big)}^\star =
{F_h\phantom{|}}^{\!\!\!\star} \, $.  Applying Proposition 2.2
to the pair  $ \, (H,K) = (U_h,F_h) \, $  we  get  $ \, {U_h}' =
{\big( {F_h}^{\!\!\vee} \big)}^* \, $  and  $ \, {F_h}^{\!\!\vee}
= {\big( {U_h}' \,\big)}^\star \, $.  By Proposition 2.4 (as  $ \,
F_h \in \QFSHA \subseteq \HA_{\otimestilde}^{w-I} \, $,  \, by
Remark 1.3{\it (c)\/})  and by Proposition 2.6 we can apply
Proposition 2.2 to the pair  $ \, (H,K) = \big( {F_h}^{\!\!\vee},
{U_h}' \,\big) \, $,  \, thus getting  $ \, {\big( {U_h}'
\,\big)}^{\!\vee} = {\Big( \big( {F_h}^{\!\!\vee} \big)'
\Big)}^{\!\star} = {F_h\phantom{|}}^{\!\!\!\star} = U_h \, $.
$ \square $
\enddemo

\vskip7pt

\proclaim{Lemma 3.8} Drinfeld's functors on quantum groups
preserve equivalences: if  $ \, H_1 \equiv H_2 \, $
in \QFSHA, resp.~in \QUEA, then  $ \, {H_1}^{\!\vee} \equiv
{H_2}^{\!\vee} \, $  in \QUEA, resp.~$ \, {H_1}' \equiv {H_2}'
\, $  in \QFSHA.
\endproclaim

\demo{Proof}  Let  $ \, H_1 $,  $ H_2 \in \QFSHA \, $  be two
equivalent quantisations of some  $ F[[\gerg]] $,  and identify
them   --- as  $ \Bbbk[[h]] $--modules  ---   with  $ \, H :=
F[[\gerg]][[h]] \, $,  \, so that the equivalence  $ \, \phi
\, \colon \, H = H_1 \equiv H_2 = H \, $  reads  $ \, \phi =
\hbox{\it id}_{\scriptscriptstyle H} + h \, \phi_+ \, $  for some
$ \, \phi_+ \in \hbox{\it End}_{\,\Bbbk[[h]]}(H) \, $.  By definition,
$ \, \phi_+ = (\phi - \hbox{\it id}_{\scriptscriptstyle H}) \big/
h \, $;  \, therefore, for all  $ \, n \in \N $,  \, we have
  $$  \big( \phi^{\otimes n} - {\hbox{\it id}_{\scriptscriptstyle
H}}^{\hskip-3pt \otimes n} \big) \big/ h = \Bigg( \sum_{k=0}^{n-1}
\phi^{\otimes k} \otimes (\phi - \hbox{\it id}_{\scriptscriptstyle H})
\otimes {\hbox{\it id}_{\scriptscriptstyle H}}^{\hskip-3pt \otimes
(n-k-1)} \Bigg) \Bigg/ h = \sum_{k=0}^{n-1} \phi^{\otimes k} \otimes
\phi_+ \otimes {\hbox{\it id}_{\scriptscriptstyle H}}^{\hskip-3pt
\otimes (n-k-1)} \, .  $$
Now, let  $ \, J := \text{\it Ker}\,(\epsilon_{\scriptscriptstyle H})
\, $:  \, since  $ \phi $  is a Hopf isomorphism, it maps  $ J $  into
itself, hence also  $ \, \phi_+(J) = \left( \! {\; \phi - \hbox{\it
id}_{\scriptscriptstyle H} \, \over h } \! \right)(J) \subseteq
J \, $.  Letting  $ \, m_n \, \colon \, H^{\otimes n} \!\rightarrow
H \, $  be the  $ n $--fold  \hbox{multiplication, we have}
  $$  \displaylines{
   \quad   \phi_+ \big( J^n \big) = \big( (\phi -
\hbox{\it id}_{\scriptscriptstyle H}) \big/ h \big)
\big( J^n \big) = m_n \Big( \big((\phi^{\otimes n} -
{\hbox{\it id}_{\scriptscriptstyle H}}^{\hskip-3pt \otimes n})
\big/ h \big) \big( J^{\otimes n} \big) \Big) =   \hfill  \cr
   \hfill   = m_n \Big( \sum_{k=0}^{n-1} \phi^{\otimes k} \otimes
\phi_+ \otimes {\hbox{\it id}_{\scriptscriptstyle H}}^{\hskip-3pt
\otimes (n-k-1)} \big( J^{\otimes n} \big) \Big) \subseteq m_n
\big( J^{\otimes n} \big) = J^n \, ,  \cr }  $$
i.e.~$ \, \phi_+ \big( J^n \big) \subseteq J^n \, $  for all  $ n $,
so  $ \, \phi_+^\vee \big( H^\vee \big) \subseteq H^\vee \, $,  \,
where  $ \, \phi_+^\vee \, $  is the extension of  $ \phi_+ $  to
$ H^\vee $.  Thus  $ \, \phi^\vee = \hbox{\it id}_{\scriptscriptstyle
H^\vee} + h \, \phi_+^\vee \, $  with  $ \, \phi_+^\vee \in \hbox{\it
End}_{\,\Bbbk[[h]]}\big( H^\vee \big) \, $,  \, so  $ \phi^\vee $  is an
equivalence in  \QUEA.
                                                 \par
   Similarly, let  $ \, H_1 $,  $ H_2 \in \QUEA \, $  be two
equivalent quantisations of some  $ U(\gerg) $,  and iden\-tify
them   --- as  $ \Bbbk[[h]] $--modules  ---   with  $ \, H :=
U(\gerg)[[h]] $.  Then the equivalence $ \, \varphi \,
\colon \, H_1 \equiv H_2 \, $  reads  $ \, \varphi =
\hbox{\it id}_{\scriptscriptstyle H} + h \, \varphi_+ \, $
for some  $ \, \varphi_+ \in \hbox{\it End}_{\,\Bbbk[[h]]}(H) \, $.
As above,  $ \, \big( \varphi^{\otimes n} - {\hbox{\it
id}_{\scriptscriptstyle H}}^{\hskip-3pt \otimes n} \big)
\big/ h = \sum_{k=0}^{n-1} \varphi^{\otimes k} \otimes \varphi_+
\otimes {\hbox{\it id}_{\scriptscriptstyle H}}^{\hskip-3pt
\otimes (n-k-1)} \, $,  \, so  $ \, \delta_n \circ \varphi
= \varphi^{\otimes n} \circ \delta_n \, $  (for  $ \varphi $
is a  {\sl Hopf\/}  isomorphism!), hence
  $$  \delta_n \circ \varphi_+ = \Big( \big( \varphi^{\otimes n}
- {\hbox{\it id}_{\scriptscriptstyle H}}^{\hskip-3pt \otimes n}
\big) \big/ h \Big) \circ \delta_n = \sum\nolimits_{k=0}^{n-1} \big(
\varphi^{\otimes k} \otimes \varphi_+ \otimes {\hbox{\it
id}_{\scriptscriptstyle H}}^{\hskip-3pt \otimes (n-k-1)} \big)
\circ \delta_n  $$
for all  $ \, n \in \N \, $.  Therefore
  $$  \displaylines{
   \quad   \delta_n \big( \varphi_+ (H') \big) = \sum_{k=0}^{n-1}
\big( \varphi^{\otimes k} \otimes \varphi_+ \otimes
{\hbox{\it id}_{\scriptscriptstyle H}}^{\hskip-3pt \otimes
(n-k-1)} \big) \big( \delta_n(H') \big) \subseteq   \hfill  \cr
   \hfill   \subseteq \sum\nolimits_{k=0}^{n-1} \big(
\varphi^{\otimes k} \otimes \varphi_+ \otimes
{\hbox{\it id}_{\scriptscriptstyle H}}^{\hskip-3pt \otimes
(n-k-1)} \big) \big( h^n H^{\otimes n} \big) \subseteq
h^n H^{\otimes n}  \cr }  $$
that is  $ \, \delta_n \big( \varphi_+ (H') \big) \subseteq h^n
H^{\otimes n} \, $  for all  $ \, n \in \N \, $,  \, so  $ \,
\varphi_+ (H') \subseteq H' \, $;  \, hence  $ \, \varphi' :=
\varphi{\big\vert}_{H'} = \hbox{\it id}_{\scriptscriptstyle H'}
+ h \, \varphi_+{\big\vert}_{H'} \, $  with  $ \,
\varphi_+{\big\vert}_{H'} \in \hbox{\it End}_{\,\Bbbk[[h]]}
\big( H' \big) \, $,  \, so that  $ \varphi' $  is an
equivalence in  $ \QFSHA $.   $ \square $
\enddemo

\vskip7pt

   Finally, our efforts are rewarded:

\vskip9pt

\demo{$ \underline{\hbox{\it Proof of Theorem 1.6}} $}  It is
enough to collect together the previous results.  Proposition 3.2
and 3.6 together with Lemma 2.8 ensure that the functors in the
claim are well-defined, and that relations  $ \circledast $  do
hold.  Proposition 3.4 and 3.7 show these functors are inverse to
each other.  Finally, Lemma 3.8 prove that they preserve equivalence.
$ \square $
\enddemo

\vskip7pt

  {\bf 3.9 Generalizations.} \, In this paper we dealt with finite
dimensional Lie bialgebras.  What about the infinite dimensional
case?  Hereafter we sketch a draft of an answer.
                                            \par
   Let  $ \gerg $  be an infinite dimensional Lie bialgebra; then
its linear dual  $ \gerg^* $  is a Lie bialgebra only in a  {\sl
topological\/}  sense: in fact, the natural Lie cobracket takes values
in the "formal tensor product"  $ \, \gerg^* \dot{\otimes} \gerg^*
:= {(\gerg \otimes \gerg)}^* $,  \, which is the completion of
$ \gerg^* \otimes \gerg^* $  w.r.t.~the weak topology.  Note that
a vector subspace  $ \gerg^\times $  of  $ \gerg^* $  is dense in
$ \gerg^* $  w.r.t.~the weak topology if and only if the restriction
$ \, \gerg \times \gerg^\times \! \rightarrow \Bbbk \, $  of the natural
evaluation pairing is perfect.
                                            \par
   If  $ \gerg $  is a Lie bialgebra in the strict algebraic
sense (i.e.~$ \delta_{\gerg} \subseteq \gerg \otimes \gerg $)
then  $ U(\gerg) $  is a co-Poisson Hopf algebra as usual; if
instead  $ \gerg $  is a Lie bialgebra in the topological sense
(i.e.~$ \delta_{\gerg} \subseteq \gerg \dot{\otimes} \gerg $)
then  $ U(\gerg) $  is a topological co-Poisson Hopf algebra,
whose co-Poisson bracket takes values in a suitable completion
$ \, U(\gerg) \,\dot{\otimes}\, U(\gerg) \, $  of  $ \, U(\gerg)
\otimes U(\gerg) \, $.  On the other hand, for any Lie bialgebra
$ \gerg $  (both algebraic or topological) we can consider  {\sl
two\/}  objects to play the role of  $ F[[\gerg]] $,  namely  $ \,
\Fstel[[\gerg]] := {U(\gerg)}^* \, $  (the linear dual of $ U(\gerg)
\, $),  endowed with the weak topology, and  $ \Finfty[[\gerg]] $,
the  $ {\frak m}_e $--adic  completion  of  $ F[G] $   --- provided
the latter exists! ---   at the maximal ideal  $ {\frak m}_e $  of
$ \, e \in G \, $,  with the  $ {\frak m}_e $--adic  topology.
Both  $ \Fstel[[\gerg]] $  and  $ \Finfty[[\gerg]] $  are  {\sl
topological\/}  Poisson Hopf algebras (the coproduct taking values
in a suitable topological tensor product), complete w.r.t.~to their
topology.  Moreover, there are natural pairings of (topological) Hopf
algebras between $ U(\gerg ) $  and  $ \Fstel[[\gerg]] $  and between
$ U(\gerg ) $  and  $ \Finfty[[\gerg]] $,  compatible with the Poisson
and co-Poisson structures.  We still have  $ \, \Fstel[[\gerg]]
\supseteq \Finfty[[\gerg]] \, $,  \, but contrary to the finite
dimensional case we may have  $ \, \Fstel[[\gerg]] \not=
\Finfty[[\gerg]] \, $.
                                            \par
   Let  $ \HA_{\otimeshat} $,  resp.~$ \HA_{\otimestilde} $,  be
defined as in \S 1.1.  In addition, define  $ \HA_{\otimesbar} $  to
be the tensor category of all (topological) Hopf  $ \Bbbk[[h]] $--algebras
$ H $  such that: \;  {\it (a)} \,  $ H $  is complete w.r.t.~the
$ I_{\scriptscriptstyle H} $--adic  topology; \;  {\it (b)}  the
tensor product  $ H_{\!1} \! \otimesbar \! H_2 $  is the completion
of the algebraic tensor product  $ H_{\!1} \otimes_{\Bbbk[[h]]} H_2 $
w.r.t.~the  $ I_{\scriptscriptstyle H_1 \otimes_{\Bbbk[[h]]} H_2} $--adic
topology; in particular, the coproduct of  $ H $  takes values in
$ \, H \otimesbar H \, $.  Then we call  $ \QUEA^\wedge $,
resp.~$ \QFSHA^\circledast $,  resp.~$ \QFSHA^\infty $,
the subcategory of  $ \HA_{\otimeshat} $,  resp.~of
$ \HA_{\otimestilde} $,  resp.~of  $ \HA_{\otimesbar} $,
composed of all objects whose specialisation at  $ \, h = 0 \, $
is isomorphic to some  $ U(\gerg) $,  resp.~some  $ \Fstel[[\gerg]] $,
resp.~some  $ \Finfty[[\gerg]] $;  here  $ \gerg $  is any Lie
bialgebra.  However, note that if  $ \, H \in \QUEA^\wedge \, $
then the Poisson cobracket  $ \delta $  of its semiclassical limit
$ \, H_0 = U(\gerg) \, $  (defined as in  Remark 1.3{\it (a)})  takes
values in  $ \, H_0 \otimes H_0 \, $,  \, so that  $ H_0 $  is an
{\sl algebraic\/}  (not topological) co-Poisson Hopf algebra hence
$ \gerg $  is an  {\sl algebraic\/}  Lie bialgebra; this means that if
we start instead from a  {\sl topological\/}  Lie bialgebra  $ \gerg $
we cannot quantize  $ U(\gerg) $  in the category \QUEA: what's wrong
is the tensor product  $ \otimeshat $  because, roughly,  $ \, H^\vee
\otimeshat H^\vee \, $  is "too small"!  Thus one must define a new
category  $ {\Cal T}_{\otimescheck} $  with the same objects than
$ {\Cal T}_{\otimeshat} $  but with a "larger" tensor product  $ \,
\otimescheck \, $  (a suitable completion of  $ \, \otimes_{\Bbbk[[h]]}
\, $)  and then consider the tensor category  $ \HA_{\otimescheck} $
of all (topological) Hopf algebras in  $ {\Cal T}_{\otimescheck} $,
and the subcategory  $ \QUEA^\vee $  whose objects have some
$ U(\gerg) $  as specialisation at  $ \, h = 0 \, $:  \, then
in this case the Lie bialgebra $ \gerg $  might be of  {\sl
topological\/}  type as well.
                                            \par
   Now let's have a look back.  We review our previous work and,
somewhat roughly, point out how far (and in which way) its results
extend to the more generals setting.

\vskip3pt

   $ \underline{\hbox{\sl Lemma 2.1}} \, $:  \, This turns into:
{\it Dualisation  $ \, H \mapsto H^* $,  resp.~$ \, H \mapsto
H^\star $,  defines a contravariant functor of tensor categories
$ \; \HA_{\otimeshat} \! \loongrightarrow \HA_{\otimestilde} \, $,
\; resp.~$ \; \HA_{\otimestilde} \! \loongrightarrow \HA_{\otimeshat}
\, $,  \; which, if  $ \text{\it char}\,(\Bbbk) \! = \! 0 $,  restrict to
$ \; \QUEA^\wedge \! \longrightarrow \QFSHA^\circledast \, $,  \;
resp.~$ \; \QFSHA^\circledast \! \longrightarrow \QUEA^\wedge \, $.}

\vskip5pt

   Indeed, this suggest to define  $ \otimescheck $  in such a
way that dualisation  $ \, H \mapsto H^* $,  resp.~$ \, H \mapsto
H^\star $,  defines a functor of tensor categories  $ \;
\HA_{\otimescheck} \! \longrightarrow \HA_{\otimesbar} \, $,  \;
resp.~$ \; \HA_{\otimesbar} \! \longrightarrow \HA_{\otimescheck}
\, $;  \; then, if  $ \text{\it char}\,(\Bbbk) \! = \! 0 $,  this will
restrict to  $ \; \QUEA^\vee \! \longrightarrow \QFSHA^\infty \, $,
\; resp.~$ \; \QFSHA^\infty \! \longrightarrow \QUEA^\vee \, $.

\vskip5pt

   $ \underline{\hbox{\it Proposition 2.2}} \, $:  \,  {\it This still
holds true for any pair  $ \, (H,K) \in \HA_{\otimeshat} \! \times
\HA_{\otimestilde} \, $.  Moreover, it should also holds true for
any pair  $ \, (H,K) \in \HA_{\otimescheck} \! \times \HA_{\otimesbar}
\, $  of Hopf algebras in duality.}

\vskip5pt

   $ \underline{\hbox{\sl Lemma 2.3}} $  \, {\it still holds
true up to replacing  $ \HA_{\otimestilde}^{w-I} $  with
$ \HA_{\otimesbar} $.}

\vskip5pt

   $ \underline{\hbox{\it Proposition 2.4}} $  \,  {\it still holds
true but for replacing  $ \HA_{\otimestilde}^{w-I} $  in part (b)
with  $ \HA_{\otimesbar} $.}

\vskip5pt

   On the other hand, if we consider  $ H^\vee $  for any  $ \, H \in
\HA_{\otimestilde} \, $  then the sole thing which goes wrong is that
$ \, \Delta\big(H^\vee\big) \not\subseteq H^\vee \otimeshat H^\vee $,
\, in general: indeed,  $ \, \Delta\big(H^\vee\big) \, $  will lie
in something larger.  Well, the definition of the category
$ \HA_{\otimescheck} $  above should fit in this frame to
give exactly  $ \, H^\vee \in \HA_{\otimescheck} \, $.  Once
one has fixed this point,  {\it our arguments still prove that
$ \, H^\vee \in \QUEA^\vee \, $};  \, thus one has a further
version of Proposition 2.4, and similarly a proper version
of Lemma 2.3 should hold with  $ \HA_{\otimestilde} $,
resp.~$ \HA_{\otimescheck} $,  instead of  $ \HA_{\otimesbar} $,
resp.~$ \HA_{\otimeshat} $.
                                            \par
   In any case, we can also drop at all the question of what kind
of Hopf algebra  $ H^\vee $  is, for in any case  {\sl the proof
of Proposition 2.4 will always prove\/}  the following:

\vskip3pt

   {\it  If  $ \, H \in \HA_{\otimestilde} \, $,  \, then  $ \,
{H^\vee}_{\!0} = U(\gerg) \, $  for some Lie bialgebra (perhaps
of topological type).}

\vskip5pt

   $ \underline{\hbox{\it Proposition 2.6}} $,
$ \underline{\hbox{\it Lemma 2.7}} \, $:  \, The
proofs we give actually show the following:
                                       \hfill\break
  {\it If  $ \, H \in \HA_{\otimeshat} $,  \,
then  $ \, H' \in \QFSHA^\circledast $.  If
$ \, H_1 $,  $ H_2 \in \HA_{\otimeshat} $,  \, then
$ \, {\big( H_1 \otimeshat H_2 \big)}' = {H_1}^{\!\prime}
\otimestilde {H_2}^{\!\prime} \, $.}

\vskip5pt

   To prove these results we used a duality argument, relying on
Lemma 2.1.  Alternatively, given  $ \, H \in \HA_{\otimeshat} \, $
or  $ \, H \in \HA_{\otimescheck} \, $  we can prove as before
that  $ H' $  is a unital  $ \Bbbk[[h]] $--subalgebra  of  $ H $,
and  {\sl also\/}  that  $ H' $  is complete w.r.t.~the  $ \,
I_{\scriptscriptstyle H'} $--adic  topology and is closed for the
antipode; then what one misses to have  $ \, H' \in \HA_{\otimesbar}
\, $  is a control on  $ \Delta\big(H'\big) \, $.  Moreover, one
proves as before that  $ \, {\big(H'\big)}_0 \, $  is commutative.

\vskip5pt

   $ \underline{\hbox{\it Lemma 2.8}} \, $:  \, The way the action
of Drinfeld's functors on morphisms is defined here still works for
any one of the categories we are considering now.

\vskip5pt

   $ \underline{\hbox{\it Sections 3.1 through 3.4}} \, $:
\, {\it These results also hold in a greater generality}.
                                                \par
   Indeed, changing a few details we can adapt the discussion in \S 3.1
and the (claim and) proof of Proposition 3.2 and of Proposition 3.4
(noting that Lemma 3.3, which still holds true untouched) to the case
of  $ \, \Fhinfty[[\gerg]] \in \QFSHA^\infty \, $.  The outcome is

\vskip2pt

   {\it If  $ \; \Fhinfty[[\gerg]] \in \QFSHA^\infty \, $,  \; then
$ \; {\Fhinfty[[\gerg]]}^\vee \! \in \QUEA^\wedge \, $,  \, namely
$ \; {\Fhinfty[[\gerg]]}^\vee \! = U_h(\gerg^\times) \, $  where
$ \gerg^\times $  is an  {\sl algebraic}  Lie bialgebra which embeds
in  $ \gerg^* $  as a dense Lie sub-bialgebra.  Moreover,  $ \, {\big(
{\Fhinfty[[\gerg]]}^\vee \big)}' = {\Fhinfty[[\gerg]]} \, $.}

\vskip3pt

   Similarly, one can apply the same arguments to  $ \,
\Fhstel[[\gerg]] \in \QFSHA^\circledast \, $  and get
essentially the same result but with  $ \QUEA^\vee $
instead of  $ \QUEA^\wedge $.  Then again the sole real
problem is to provide a proper definition for the category
$ \HA_{\otimescheck} $  (or at least  $ \QUEA^\vee $)  which
fit well with these results.  Once this (non-trivial...) point
is set, the result would read

\vskip2pt

   {\it If  $ \, \Fhstel[[\gerg]] \! \in \! \QFSHA^\circledast $,
\, then  $ \, {\Fhstel[[\gerg]]}^\vee \in \QUEA^\wedge \, $,  \,
namely  $ \, {\Fhstel[[\gerg]]}^\vee \!\! = U_h(\gerg^*) $,  \,
where  $ \gerg^* $  is the dual  {\sl (topological)}  Lie bialgebra
to  $ \gerg $.  Moreover,  $ \, {\big( {\Fhstel[[\gerg]]}^\vee \big)}'
= {\Fhstel[[\gerg]]} \, $.}

\vskip5pt

   $ \underline{\hbox{\it Sections 3.5 through 3.7}} \, $:  \, {\it
These again hold in a greater generality}.
                                                \par
   In this case, the main tool is the use duality functor to switch
from QUEA to QFSHA and the property of Drinfeld's functors of being
dual to each other ensured by Proposition 2.2.  Therefore, our
arguments apply  {\it verbatim\/}  to the case of  $ \, U_h(\gerg)
\in \QUEA^\wedge $.  As for  $ \, U_h(\gerg) \in \QUEA^\vee $,  \,
everything goes true as well the same  {\sl provided\/}  Lemma 2.1
and Proposition 2.2 have been properly extended to deal with
$ \QUEA^\vee $  and  $ \QFSHA^\infty $,  as mentioned above.

\vskip5pt

   $ \underline{\hbox{\it Lemma 3.8}} \, $:  \, Here again (as for
Lemma 2.8) our analysis still works for any one of the categories
we are considering now.

\vskip5pt

   In a nutshell, we can say that, up to some details to be fixed,

\vskip3pt

   {\it The quantum duality principle holds, in a suitable
formulation, also for infinite dimensional Lie bialgebras,
both algebraic and topological.}

\vskip7pt

   {\bf 3.10 Examples.} \, Several examples about finite dimensional
Lie bialgebras can be found in [Ga2]: there we consider quantum groups
"\`a la Jimbo-Lusztig", but one can easily translate all definitions
and results into the language
              "\`a la Drinfeld'" we use in the present paper.\break
\indent   We consider now some infinite dimensional samples.  Let
$ \gerg $  be a simple finite dimensional complex Lie algebra, and
$ \hat{\gerg} $  the associated untwisted affine Kac-Moody algebra,
with the well-known Sklyanin-Drinfeld structure of Lie bialgebra; let
also  $ \hat{\gerh} $  be defined as in [Ga1], \S 1.2.  Then both
$ \hat{\gerg}^* $ and  $ \hat{\gerh} $  are  {\sl topological\/}  Lie
bialgebras, with  $ \hat{\gerh} $  dense inside  $ \hat{\gerg}^* $.
                                            \par
  Consider the quantum groups  $ {\Cal U}^{\scriptscriptstyle M}
(\hat{\gerg}) $  and  $ {\frak U}^{\scriptscriptstyle M}
(\hat{\gerg}) $  defined in [Ga1], \S 3.3.  We can reformulate
the definition of the first in Drinfeld's terms via the usual
"dictionary": pick generators  $ \, H_i = \log \big( K_i \big)
\big/ \log(q_i) \, $  instead of the  $ K_i $'s,  take  $ \, h
= \log(q) \, $  and fix  $ \Bbbk[[h]] $  as ground ring, and finally
complete the resulting algebra w.r.t.~the  $ h $--adic  topology;
then we have exactly  $ \, {\frak U}^{\scriptscriptstyle M}
(\hat{\gerg}) \in \HA_{\otimeshat} \, $  and  $ \,
{\frak U}^{\scriptscriptstyle M}(\hat{\gerg}) \,{\buildrel
{h \rightarrow 0} \over \llongrightarrow}\; U(\hat{\gerg}) \, $,
\, so  $ \, U_h(\hat{\gerg}) := {\frak U}^{\scriptscriptstyle M}
(\hat{\gerg}) \in \QUEA^\wedge \, $  (discarding the choice of
the weight lattice  $ M $).  On the other hand, doing the same
"translations" for  $ {\Cal U}^{\scriptscriptstyle M}(\hat{\gerg}) $
and completing w.r.t.~the weak topology  {\sl or\/}  w.r.t.~the
$ I $--adic  topology we obtain  {\sl two\/}  different objects, in
$ \HA_{\otimestilde} $  and in  $ \HA_{\otimesbar} $  respectively,
with semiclassical limit  $ \Fstel[[\hat{\gerg}^*]] $  and
$ \Finfty\big[\big[\hat{\gerh}\big]\big] $  respectively; then
we call them  $ \Fhstel[[\hat{\gerg}^*]] $  and  $ \Fhinfty
\big[\big[\hat{\gerh}\big]\big] $  respectively, with  $ \,
\Fhstel[[\hat{\gerg}^*]] \in \QFSHA^\circledast \, $  and
$ \, \Fhinfty\big[\big[\hat{\gerh}\big]\big] \in \QFSHA^\infty \, $.
Now, acting as outlined in [Ga2], \S 3, one finds  $ \, {\Fhinfty
\big[\big[\hat{\gerh}\big]\big]}^\vee \! = U_h(\hat{\gerg}) \, $
and  $ \, {U_h(\hat{\gerg})}' = \Fhinfty\big[\big[\hat{\gerh}\big]\big]
\, $,  whilst  $ \, {\Fhstel[[\hat{\gerg}^*]]}^\vee \, $  instead is
a suitable completion of  $ U_h(\hat{\gerg}) $,  which should be
an object of  $ \QUEA^\vee $:  indeed, we have  $ \, {\Fhstel
[[\hat{\gerg}^*]]}^\vee {\buildrel {h \rightarrow 0} \over
\llongrightarrow}\; U\big(\hat{\gerh}^*\big) \, $,  \,
and  $ \hat{\gerh}^* $  is a topological Lie bialgebra
in perfect duality with  $ \hat{\gerg}^* $.
                                            \par
  Dually, consider the quantum groups  $ {\Cal U}^{\scriptscriptstyle
M} \big(\hat{\gerh}\big) $  and  $ {\frak U}^{\scriptscriptstyle M}
\big(\hat{\gerh}\big) $  defined in [Ga1], \S 5; as above we
can rephrase their definition, and then we find the following.
First, the formul\ae{} for the coproduct imply that  $ \,
{\Cal U}^{\scriptscriptstyle M}\big(\hat{\gerh}\big) \not\in
\HA_{\otimesbar} \, $  {\sl but}  $ \, {\Cal U}^{\scriptscriptstyle M}
\big(\hat{\gerh}\big) \in \HA_{\otimestilde} \, $,  \, and  $ \,
{\Cal U}^{\scriptscriptstyle M}\big(\hat{\gerh}\big) \,{\buildrel
{h \rightarrow 0} \over \llongrightarrow}\; \Fstel \big[\big[
\hat{\gerh}^* \big]\big] \, $,  \, thus  $ \, \Fhstel \big[\big[
\hat{\gerh}^* \big]\big] := {\Cal U}^{\scriptscriptstyle M}
\big(\hat{\gerh}\big) \in \QFSHA^\circledast \, $.  Second,
$ \, {\frak U}^{\scriptscriptstyle M}\big(\hat{\gerh}\big)
\not\in \HA_{\otimeshat} \, $  {\sl but}  $ \,
{\frak U}^{\scriptscriptstyle M}\big(\hat{\gerh}\big)
\in \HA_{\otimescheck} \, $,  \, and  $ \,
{\frak U}^{\scriptscriptstyle M}\big(\hat{\gerh}\big)
\,{\buildrel {h \rightarrow 0} \over \llongrightarrow}\;
U\big(\hat{\gerh}\big) \, $,  \, so in fact  $ \,
U_h\big(\hat{\gerh}\big) := {\frak U}^{\scriptscriptstyle M}
\big(\hat{\gerh}\big) \in \QUEA^\vee \, $.  By an analysis
like that in  [Ga2]  one shows also that  $ {\Fhstel\big[\big[
\hat{\gerh}^* \big]\big]}^\vee $  is an object of  $ \QUEA^\vee $,
it is a suitable completion of  $ U_h\big(\hat{\gerh}\big) $,
and  $ \, {\Fhstel \big[\big[\hat{\gerh}^*\big]\big]}^\vee
{\buildrel {h \rightarrow 0} \over \llongrightarrow}\;
U \big( \hat{\gerg}^* \big) \, $;  \, moreover,  $ \,
{\Big( {\Fhstel \big[\big[\hat{\gerh}^*\big]\big]}^\vee \Big)}'
= \Fhstel\big[\big[\hat{\gerh}^*\big]\big] \, $.  On the other
hand, one has  $ \, {U_h\big(\hat{\gerh}\big)}' = \Fhstel\big[\big[
\hat{\gerh}^* \big]\big]  \in \QFSHA^\circledast $,  \, and so
$ \, {\Big({U_h\big(\hat{\gerh}\big)}'\Big)}^\vee = U_h \big(
\hat{\gerh} \big) \, $.

\vskip1,1truecm

\vskip0,3truecm


\Refs
  \widestnumber\key {FRT}

\vskip3pt

\ref
 \key  A   \by  N. Abe
 \book  Hopf algebras
 \publ  Cambridge Tracts in Mathematics  {\bf 74}
 \publaddr  Cambridge University Press, Cambridge   \yr  1980
\endref

\ref
 \key  CG   \by  N. Ciccoli, F. Gavarini
 \paper  A quantum duality principle for Poisson homogeneous spaces
 \jour  Preprint
\endref

\ref
 \key  Dr   \by  V. G. Drinfeld
 \paper  Quantum groups
 \inbook  Proc. Intern. Congress of Math. (Berkeley, 1986)  \yr  1987
 \pages  798--820
\endref

\ref
 \key  E   \by  B. Enriquez
 \paper  Quantization of Lie bialgebras and shuffle algebras
of Lie algebras
 \jour  Preprint  \break  math.QA/0008128
 \yr  2000
\endref

\ref
 \key  EK1   \by  P. Etingof, D. Kazhdan
 \paper  Quantization of Lie bialgebras, I
 \jour  Selecta Math. (New Series)   \vol  2   \yr  1996   \pages  1--41
\endref

\ref
 \key  EK2   \by  P. Etingof, D. Kazhdan
 \paper  Quantization of Poisson algebraic groups and Poisson
homogeneous spaces
 \book  Sym\'etries quantiques (Les Houches, 1995)
 \publ  North-Holland   \publaddr  Amsterdam
 \yr  1998   \pages  935--946
\endref

\ref
 \key  FRT   \by  L. D. Faddeev, N. Yu. Reshetikhin, L. A. Takhtajan
 \paper  Quantum groups
 \jour  in: M. Kashiwara, T. Kawai (eds.),  {\it Algebraic Analysis},
 \publ  Academic Press   \publaddr  Boston
 \yr 1989  \pages 129--139
\endref

\ref
 \key  KT   \by  C. Kassel, V. Turaev
 \paper  Biquantization of Lie bialgebras
 \jour  Pac. Jour. Math.   \vol  195   \yr  2000   \pages  297--369
\endref

\ref
 \key  Ga1   \by  F. Gavarini
 \paper  Dual affine quantum groups
 \jour  Math. Zeitschrift   \vol  234   \yr  2000   \pages  9--52
\endref

\ref
 \key  Ga2   \by  F. Gavarini
 \paper  The global quantum duality principle: theory,
examples, applications
 \jour  Preprint  \break  math.QA/0108015   \yr  2001
\endref

\ref
 \key  M   \by  S. Montgomery
 \book  Hopf Algebras and Their Actions on Rings
 \publ  CBMS Regional Conference Series in Mathematics  {\bf 82},
American Mathematical Society
 \publaddr  Providence, RI   \yr  1993
\endref


\endRefs
\enddocument